\documentclass{amsart}
\usepackage{amssymb}

\newtheorem{theorem}{Theorem}
\newtheorem{proposition}[theorem]{Proposition}
\newtheorem{definition}[theorem]{Definition}
\newtheorem{lemma}[theorem]{Lemma}
\newtheorem{corollary}[theorem]{Corollary}

\numberwithin{equation}{section}
\numberwithin{theorem}{section}

\begin{document}
\title[$L^p$ Modules]{NONCOMMUTATIVE $L^p$ MODULES}
\author{Marius Junge}
\author{David Sherman}
\address{Department of Mathematics\\ University of Illinois\\ Urbana, IL 61801}
\email[Marius Junge]{junge@math.uiuc.edu}
\email[David Sherman]{dasherma@math.uiuc.edu}
\begin{abstract}
We construct classes of von Neumann algebra modules by considering ``column
sums" of noncommutative $L^p$ spaces. ÊOur abstract characterization is based on
an $L^{p/2}$-valued inner product, thereby generalizing Hilbert C*-modules and
representations on Hilbert space. ÊWhile the (single) representation theory is
similar to the $L^2$ case, the concept of $L^p$ bimodule ($p \ne 2$) turns out to be nearly
trivial.
\end{abstract}

\maketitle

\section{Introduction}

Noncommutative $L^p$ spaces, by now, are standard objects in the theory of
operator algebras. ÊStarting with a von Neumann algebra $\mathcal{M}$, there are
a variety of equivalent methods for producing the (quasi-)Banach space
$L^p(\mathcal{M}).$ ÊIf $\mathcal{M}$ is $L^\infty(X, \mu)$, the result is
(isometric to) $L^p(X, \mu)$, so this can rightfully be thought of as a
generalization to noncommutative measure spaces. ÊWhen $\mathcal{M}$ is
semifinite, the presence of a trace offers great simplification, but in general
one needs modular theory [H].

These spaces have many aspects worthy of investigation. ÊAs Banach spaces, their isometries have been investigated by many authors [Ye], [W2]; others have used the matrix order [Sc] or operator space techniques [JNRX]. Ê(For a more complete bibliography see [PX].) ÊWe focus here on the
\textit{module} structure. ÊIndeed, the inclusion as left (or right) multipliers
$$\mathcal{M} \hookrightarrow \mathcal{B}(L^p(\mathcal{M}))$$
is isometric. ÊIf Hilbert space representations are (categorically) generated by
$L^2(\mathcal{M})$, and self-dual C*-modules are generated by
$L^\infty(\mathcal{M}) = \mathcal{M}$, where are the modules generated by
$L^p(\mathcal{M})$? ÊThis paper sets out to describe the missing $L^p$
representation theory.

Proceeding by analogy, our target is the class of ``columns of $L^p(\mathcal{M})$". ÊWe show that a sufficient condition for an $\mathcal{M}$-module to
belong to this class is the existence of an
$L^{p/2}(\mathcal{M})$-\textit{valued inner product}. ÊThe description which
results is a natural generalization of the cases $p=2$ (the usual decomposition
for Hilbert space representations) and $p=\infty$ (see [Pa]). ÊWe employ a
variety of methods, but perhaps the most notable direction is a consistent
translation of Connes' $L^2$ spatial theory [C] to the $L^p$ setting.

Building on results about the module structure of $L^p(\mathcal{M})$ which are
interesting in their own right, we find that the $L^p$ representation theory is
largely analogous to the $L^2$ case, with a well-behaved sum and relative tensor
product. ÊIt would therefore seem natural that there be a similarly rich
bimodule category, i.e. a theory of $L^p$ correspondences. ÊBut surprisingly, the category is nearly trivial: when $p \ne 2$, there is an $L^p$ $\mathcal{M}$-$\mathcal{N}$ bimodule if and only if $\mathcal{M}$ and $\mathcal{N}$ are Morita equivalent. ÊModulo a possible degeneracy where both algebras are abelian, such bimodules naturally implement an equivalence of appropriate representation categories.

Only one application - to ultraproducts - of our theory is given. ÊWe plan to
discuss further examples and development in future articles.

\section{The module structure of $L^p(\mathcal{M})$}

Throughout, $\mathcal{M}$, $\mathcal{N}$, etc. are von Neumann algebras; we frequently abbreviate
$L^p(\mathcal{M})$ to $L^p$ and understand $L^\infty(\mathcal{M})$ as
$\mathcal{M}$. ÊAll weights are normal and semifinite, so we omit the adjectives for brevity. ÊUnsubscripted $\mathfrak{H}$ denotes the separable
infinite-dimensional Hilbert space, $s(\varphi)$ is the support of
$\varphi$, and $s_\ell(x)$ (resp. $s_r(x)$) stands for the left (resp. right) support of $x$. ÊSubscripts are occasionally used to represent an action: e.g. $\mathfrak{X}_\mathcal{M}$ indicates that $\mathfrak{X}$ is a right $\mathcal{M}$-module. ÊBut when the expressions are longer, we signify a bimodule by writing out the triple: an $\mathcal{M}$-$\mathcal{N}$ bimodule $\mathfrak{X}$ is $\mathcal{M}$-$\mathfrak{X}$-$\mathcal{N}$. ÊThe phrase ``left (resp. right) action of" is freqently abbreviated to $L$ (resp. $R$) for operators or entire algebras, so that we speak of $L(x)$ or $R(\mathcal{M})$. ÊFinally, we often write $M_\infty$ for $\mathcal{B}(\mathfrak{H})$ and $M_\infty(\mathcal{M})$ for $\mathcal{B}(\mathfrak{H}) \otimes \mathcal{M}$. ÊNote that in contrast to much of the literature, the results of this paper (except for Section 6) do not require that algebras be $\sigma$-finite or that $p \ge 1$.

We assume that the reader has some basic familiarity with noncommutative $L^p$
spaces. ÊConceptually, one can think
$$L^p_+ = (L^1_+)^{1/p},$$
and we take this as a basis for our notation: the typical positive element is
$\varphi^{1/p}$, where $\varphi$ is a positive linear functional on
$\mathcal{M}$. ÊWhat this \textit{means} is a matter of perspective, as there
are many equivalent constructions of $L^p$, but we find the Haagerup
construction [H] most useful. ÊIn this setting $L^p$ is exactly the set of
$\tau$-measurable operators affiliated with the core $\widetilde{\mathcal{M}}
\simeq \mathcal{M} \rtimes_\sigma \mathbb{R}$ which are $1/p$-scaled by the dual
action: $\theta_s(T) = e^{-s/p} T$. Ê(The operator we call $\varphi$ is more
commonly called $h_\varphi.$ ÊAn unbounded weight corresponds to a positive operator satisfying all the above conditions except for $\tau$-measurability.) ÊOperator concepts like composition, positivity, left and right support,
adjoint, and polar decomposition transfer directly into the $L^p$ setting.
Basic exposition can be found in [Te], and the reader is also referred to the
elegant ``coordinate-free" approaches in [Y] (more algebraic) and [FT] (more
analytic). ÊWe use the Haagerup notation Tr for the evaluation functional on
$L^1$:
$$\text{Tr}(\omega) = \omega(1),$$
and recall that Tr implements the ``tracial" duality between $L^p$ and $L^q$:
$$ <\xi, \eta> = \text{Tr}(\xi \eta) = Ê\text{Tr}(\eta \xi), \qquad \xi \in L^p,
\eta \in L^q, \:\frac{1}{p} + \frac{1}{q} = 1.$$

In this notation,
$$\psi^{it}\varphi^{-it} = (D\psi: D\varphi)_t; \quad \varphi^{it}x\varphi^{-it}
= \sigma_t^{\varphi}(x), \qquad (t \in \mathbb{R})$$
whenever $s(\varphi)$ dominates $s(\psi), s_\ell(x), s_r(x).$ ÊThe cocycles or modular
automorphism groups extend off the imaginary line exactly when the corresponding
operator compositions do. ÊFor more discussion of negative powers of states, see [S2].

\bigskip

A fundamental fact for us is Kosaki's generalized H\"{o}lder inequality [K2]:
$$\|\xi \eta \|_r \le \|\xi\|_p \|\eta\|_q, \quad \xi \in L^p, \: \eta \in L^q,
\: \frac{1}{p} + \frac{1}{q} = \frac{1}{r}.$$
In particular, left or right multiplication by an element $x \in \mathcal{M}$ is
bounded with norm $\le \|x\|.$ ÊWe will show a stronger fact momentarily, but
first recall
\begin{equation} \label{E:order}
\varphi^{1/p} \ge \psi^{1/p} \iff \varphi^{-\frac{1}{2p}} \psi^{\frac{1}{2p}}
\text{ is a contraction in } s(\varphi)\mathcal{M}s(\varphi) \subset
\mathcal{M}.
\end{equation}
Then
\begin{equation} \label{E:ordnorm}
\varphi^{1/p} \ge \psi^{1/p} \Rightarrow \|\psi^{1/p}\| = \|(\psi^{\frac{1}{2p}}
\varphi^{-\frac{1}{2p}}) \varphi^{1/p} (\varphi^{-\frac{1}{2p}}
\psi^{\frac{1}{2p}}) \| \le \| \varphi^{1/p}\|. \end{equation}

\begin{lemma} \label{T:spec}
Let $x \ge 0,$ and consider the map $\xi \mapsto x \xi$ on $L^p$. ÊWe have
$$\inf \text{sp} (x) = \inf_{\|\xi\| = 1} \|x \xi \|, \quad \| x \|= \sup
\text{sp} (x) = \sup_{\|\xi\| = 1} \|x \xi \|.$$
\end{lemma}
\begin{proof}
We discuss $p<\infty$; $p=\infty$ only requires different wording.

H\"{o}lder's inequality is half of the last equation. ÊTo see the opposite
inequality, choose $\varepsilon$ and let $\varphi$ be a state supported on the
spectral projection $q = e([\|x\|-\varepsilon, \|x\|])$ of $x$. ÊWe use
$\eqref{E:ordnorm}$ above to get
$$\|x \varphi^{1/p} \|_p = \| \varphi^{1/p} q x^2 q \varphi^{1/p}
\|_{p/2}^{1/2}
\ge \|\varphi^{1/p} q (\|x\| - \epsilon)^2 q \varphi^{1/p} \|_{p/2}^{1/2} = (\|x\| - \varepsilon) \|\varphi^{1/p}\|_p.$$

The first equation is proven similarly.
\end{proof}
For an arbitrary element with polar decomposition $x = v|x|,$ we have $\|x \xi
\| = \| |x| \xi \|$, and the proposition alters naturally by considering the
spectrum of $|x|$.

This does not give us an ``$L^p$ spatial spectral theorem". ÊA positive operator
$x$ generates a projection-valued decomposition of the identity, and the action
on $L^p$ is still ``multiplication'' (in an appropriate sense) by $\int \lambda
de(\lambda).$ ÊBut for disjoint sets $I$ and $J$, there is no simple norm relation between the $L^p$ elements $e(I) \xi$, $e(J) \xi$ and their sum unless $p=2$. ÊThis
prevents us from using vectors to provide ($p$th roots of) measures, and we
cannot write, say, $\|\xi\|^p = \int d\|e(\lambda) \xi\|^p.$

Now we turn to a full description of the intertwiner set
$\text{Hom}(L^p_\mathcal{M}, L^q_\mathcal{M})$. ÊThe next three lemmas facilitate the proofs; the second is a slight improvement of [J], Lemma 2.3.

\begin{lemma} \label{T:dense}
For $\varphi \in \mathcal{M}_*^+$,
$$\overline{\mathcal{M} \varphi^{1/p}} =L^p s(\varphi), \qquad \overline{\varphi^{1/p} \mathcal{M}} =s(\varphi) L^p.$$
\end{lemma}

\begin{proof}
This result is well-known for $p \ge 1$, but we present a full proof for completeness.

Let $p>1$. ÊSuppose there is $\xi \in s(\varphi) L^p \setminus \overline{\varphi^{1/p} \mathcal{M}}.$ ÊBy Hahn-Banach separation we may find $\eta \in L^q$ ($p,q$ conjugate exponents) with
$$\text{Tr}(\eta \xi) >0, \qquad 0 = \text{Tr} (\eta \varphi^{1/p} y ), \quad \forall y \in \mathcal{M}.$$
Then we must have $\eta \varphi^{1/p} = 0$, so $\eta s(\varphi) = 0$. ÊBut
$$0 < \text{Tr}(\eta \xi) = \text{Tr}(\eta s(\varphi) \xi) = 0,$$
a contradiction. ÊBy a symmetric argument we have $\overline{\mathcal{M} \varphi^{1/p}} =L^p s(\varphi)$.

Keep the same $p$, and assume that we have
$$\overline{\mathcal{M} \varphi^{n/p}} =L^{p/n} s(\varphi), \qquad \overline{\varphi^{n/p} \mathcal{M}} =s(\varphi) L^{p/n}$$
for a positive integer $n$. ÊWe compute
\begin{align*}
\overline{\varphi^{(n+1)/p} \mathcal{M}} &= \overline{ \varphi^{1/p} \:\overline{ \varphi^{n/p}\mathcal{M}}} = \overline{ \varphi^{1/p} L^{p/n} }
Ê= \overline{\varphi^{1/p} \mathcal{M} L^{p/n}} \\ &= \overline{ s(\varphi) L^p \cdot L^{p/n}} = s(\varphi) L^{p/(n+1)},
\end{align*}
where the first equality is justified by H\"{o}lder: if $\varphi^{n/p}x_j$ converges, so does $\varphi^{(n+1)/p} x_j.$ ÊThe other equality is obtained similarly.

Since any positive number can be written as $p/n$ with $p>1$ and $n$ a positive integer, the result follows by induction.
\end{proof}

\begin{lemma} \label{T:strong}
Let $\varphi \in \mathcal{M}_*^+$, $p>0$, and $\{x_\alpha\}$ be a bounded net in $\mathcal{M}$. ÊIf
$$(*) \quad x_\alpha \varphi^{1/p} \to 0,$$
then $x_\alpha \xi \to 0$ for any $\xi$ in $s(\varphi)L^q$, where $q$ is any positive real. ÊThis implies that
on bounded sets, the strong topology that $\mathcal{M}$ acquires from its action
on $L^p$ does not depend on $p$.
\end{lemma}
\begin{proof}
Suppose
$$(*) \quad x_\alpha \varphi^\beta \to 0$$
for $\beta = 1/p$. ÊThen ($*$) holds for $\beta> 1/p$, as
$$ \| x_\alpha \varphi^\beta \| \le \| x_\alpha \varphi^{1/p} \|
\|\varphi^{\beta - 1/p} \| \to 0.$$
We may also conclude that ($*$) holds for $\beta = \frac{1}{2p}$ by
$$\| x_\alpha \varphi^{\frac{1}{2p}} \|^2 = \|x_\alpha \varphi^{1/p} x_\alpha^*
\| \le \| x_\alpha \varphi^{1/p} \| \|x_\alpha^* \| \to 0,$$
since $\|x_\alpha^* \|$ is bounded. ÊTogether these two steps imply that ($*$)
holds for all positive $\beta$.

Now suppose that $\xi \in s(\varphi)L^q = \overline{\varphi^{1/q}\mathcal{M}}.$ ÊGiven $\varepsilon >0$, choose $y\in \mathcal{M}$ so that
$$ \|\xi - \varphi^{1/q} y\| \le \frac{\varepsilon}{2 \sup \|x_\alpha\|}.$$
Then
\begin{align*}
\|x_\alpha \xi \| &\le \|x_\alpha \xi - x_\alpha \varphi^{1/q} y\| + \|
x_\alpha \varphi^{1/q} y \| \\ &\le \left( \sup \|x_\alpha \| \right) \| \xi -
\varphi^{1/q} y \| + \| x_\alpha \varphi^{1/q} \| \|y\|,
\end{align*}
which is less than $\varepsilon$ when $\alpha$ is so large that $\|x_\alpha
\varphi^{1/q} \| < \frac{\varepsilon}{2 \|y\|}$.

When $\mathcal{M}$ is $\sigma$-finite, this last step is the $L^p$ version of the well-known fact that for a faithful state $\varphi$,
$$ x \mapsto \varphi(x^* x)^{1/2}$$
implements the strong topology on bounded sets of $\mathcal{M}$.
\end{proof}

\begin{lemma}\label{T:addcols}
Let $\{\xi_\alpha\}\subset L^r \quad (r < \infty), \quad \{p_\alpha\} \subset
\mathcal{P}(\mathcal{M})$ be nets such that
\begin{equation}\label{E:addrows}
\xi_\alpha = \xi_\alpha p_\alpha, \quad \xi_\beta p_\alpha = \xi_\alpha \text{ for Ê}\alpha <\beta, \quad \sup_\alpha
\|\xi_\alpha\| = C < \infty.
\end{equation}
Then $\xi_\alpha$ converges in norm, say to $\xi$, and $\xi_\beta = \xi p_\beta.$
\end{lemma}

The idea is that adding columns (=increasing the right support) without
exceeding an $L^r$ bound implies convergence in $L^r$.

\begin{proof}
First we handle the case where $r > 2.$ ÊWe have $\xi_\alpha \xi_\alpha^*$ increasing and norm-bounded; let $\varphi^{2/r}$ be the weak-* limit in the reflexive Banach space
$L^{r/2}$ and write
$$\xi_\alpha \xi_\alpha^* = \varphi^{1/r} x_\alpha \varphi^{1/r} \text{ with } x_\alpha \le
q=s(\varphi).$$
Using $L^r(q\mathcal{M}q) = qL^r(\mathcal{M})q,$ weak convergence implies that
$$<\varphi^{1/r} x_\alpha \varphi^{1/r}, \psi^{1/s}> \to <\varphi^{2/r}, \psi^{1/s}>,
\qquad \forall \psi \in qL^s q \quad (\frac{2}{r} + \frac{1}{s}=1),$$
or
$$<x_\alpha, \varphi^{1/r} \psi^{1/s} \varphi^{1/r}> \to <q, \varphi^{1/r} \psi^{1/s}
\varphi^{1/r}>;$$
that is, $x_\alpha =qx_\alpha q \nearrow q$ weakly in $\mathcal{M}$. ÊNow
$$q \ge x_\alpha^{1/2} \ge x_\alpha \nearrow q \text{ weakly } \Rightarrow x_\alpha^{1/2}
\nearrow q \text{ weakly }$$
$$\Rightarrow (q-x_\alpha^{1/2})^2 = q + x_\alpha - 2x_\alpha^{1/2} \searrow 0 \text{
weakly,}$$
so $x_\alpha^{1/2} \nearrow q$ strongly. ÊBy the preceding lemma,
$$ \varphi^{1/r} x_\alpha^{1/2} \to \varphi^{1/r},$$
and therefore
$$\xi_\alpha \xi_\alpha^* = (\varphi^{1/r} x_\alpha^{1/2}) (\varphi^{1/r} x_\alpha^{1/2})^* \to
\varphi^{2/r}.$$
Finally, for $\alpha <\beta$ the increasing right supports imply
$$\|\xi_\alpha - \xi_\beta \|^2 = \|\xi_\alpha \xi_\alpha^* + \xi_\beta \xi_\beta^* - \xi_\alpha \xi_\beta^* -\xi_\beta
\xi_\alpha^*\|$$ $$ = \|\xi_\alpha \xi_\alpha^* + \xi_\beta \xi_\beta^* - \xi_\alpha \xi_\alpha^* - \xi_\alpha
\xi_\alpha^*\| = \|\xi_\beta \xi_\beta^* - \xi_\alpha \xi_\alpha^*\| \to 0.$$

If $r\le2$, still $\xi_\alpha \xi_\alpha^*$ is increasing in $L^{r/2}$ and bounded.
Choose $\gamma > 2/r$; $(\xi_\alpha \xi_\alpha^*)^{1/\gamma}$ is then norm-bounded and
increasing in a reflexive Banach space. ÊBy the above argument it converges in
norm, so the continuity of exponentiation (see [R], Lemma 3.2) implies $\xi_\alpha \xi_\alpha^* = ((\xi_\alpha \xi_\alpha^*)^{1/\gamma})^\gamma$ converges in $L^{r/2}$. ÊThe last computation of the previous paragraph again shows the convergence of $\xi_\alpha$.

Finally, set $\xi = \lim_\alpha \xi_\alpha$ and use that right multiplication by $p_\beta$ is continuous:
$$\xi p_\beta = (\lim_\alpha \xi_\alpha) p_\beta = \lim_\alpha (\xi_\alpha p_\beta) = \lim_\alpha \xi_\beta = \xi_\beta.$$
\end{proof}

When $p=\infty$, Lemma \ref{T:addcols} still holds. ÊThe same line of argument works, but instead of reflexivity one uses that von Neumann algebras are monotone closed.

The next theorem extends work of several authors and solves a problem stated in
Yamagami [Y].

\begin{theorem} \label{T:inter}
If $\frac{1}{p} + \frac{1}{r} = \frac{1}{q},$ then any bounded map in
$\text{Hom}(L^p_\mathcal{M}, L^q_\mathcal{M})$ is left composition with some
element of $L^r$.
\end{theorem}

\begin{proof}
Let $T$ be such a map. ÊIf $p=\infty$, this is easy: $T(x) = T(1)x.$ ÊSo assume
$p < \infty$, and for the moment assume $\mathcal{M}$ is $\sigma$-finite. ÊChoose a faithful $\varphi \in \mathcal{M}_*^+.$ ÊWith
$$T(\varphi^{1/p}) = v \psi^{1/q}$$
the polar decomposition, set
$$\rho^{2/p} = \varphi^{2/p} + \psi^{2/p}.$$
and write
$$\psi^{1/p} = y_1 \rho^{1/p}; \qquad \varphi^{1/p} = y_2 \rho^{1/p}$$
with $y_1, y_2$ contractive. ÊThe module property means that for any $x \in
\mathcal{M}$,
$$T(y_2 \rho^{1/p} x) = v \psi^{1/q} x = v \psi^{1/r} y_1 \rho^{1/p} x.$$
By continuity of $T$ we may conclude
$$T(y_2 \xi) = v \psi^{1/r} y_1 \xi$$
for all $\xi \in L^p$.

Now let $y_2 = |y_2^*|u$ be the polar decomposition and $q_n$ be the spectral
projection of $|y_2^*|$ corresponding to $[\frac{1}{n}, 1]$. ÊSince
$$q_n = y_2 u^* |y_2^*|^{-1} q_n,$$
$$T(q_n \xi) = T(y_2 u^* |y_2^*|^{-1} q_n \xi) = (v \psi^{1/r} y_1 u^* |y_2^*|^{-1}
q_n) (q_n \xi).$$
It follows from this that
$$\|v \psi^{1/r} y_1 u^* |y_2^*|^{-1} q_n\|_r \le \|T\|$$
for all $n$, and notice the $q_n$ are increasing to 1 since $y_2$ is
nonsingular.

If $r<\infty$, Lemma \ref{T:addcols} allows us to conclude the convergence of
this sequence; say
$$v \psi^{1/r} y_1 u^* |y_2^*|^{-1} q_n \to \eta.$$
Since $T$ agrees with $L(\eta)$ on the dense set $\cup q_n
L^p$, they are identical.

If $r = \infty$, then $\psi^{1/r}$ can be replaced with 1. ÊThe uniform bound
implies that $v y_1 u^* |y_2^*|^{-1} q_n$ converges strongly to an operator $z$
with $\|z\| \le \|T\|.$ ÊAgain, $T$ and $L(z)$ agree on $\cup q_n L^p$, so they
are identical.

\bigskip

Now we remove the $\sigma$-finiteness assumption.

Let $r<\infty.$ ÊIf $s$ is a $\sigma$-finite projection in $\mathcal{M}$, we may find a state $\varphi$ with $s(\varphi)=s$ and apply the same argument to conclude
$$T \mid_{sL^p} = L(\eta_s).$$
Then the $\eta_s$ satisfy
$$\eta_s = \eta_s s, \qquad \eta_t s = \eta_s \text{ for }s<t, \qquad \|\eta_s\| \le \|T\|.$$
Lemma \ref{T:addcols} tells us that $\eta_s$ converges along the naturally-ordered net of $\sigma$-finite projections, say to $\eta$, and $\eta_s = \eta s$. ÊFinally, if $\xi \in L^p$, $f = s_\ell(\xi)$ must be $\sigma$-finite, and
$$T(\xi) = T(f\xi) = \eta_f \xi = \eta f \xi = \eta \xi.$$

In case $r= \infty,$ the vectors $\eta_s, \eta$ are replaced by operators $z_s, z.$
\end{proof}

We single out the case $r= \infty$ as a separate corollary. ÊThough basic, there does not seem to be a proof for general $p$ in the literature. Ê(Terp [Te] settled the
case $p \ge 1$ by different methods.)

\begin{corollary} \label{T:comm}
The left and right actions of $\mathcal{M}$ on $L^p$ are commutants of each
other.
\end{corollary}

Notice that for $p \ge 1$, $(L^p)^* = L^{p'}$ can be identified with
$\text{Hom}(L_\mathcal{M}^p,L_\mathcal{M}^1)$, with $Tr$ implementing the
duality as usual. ÊIt is known [W1] that $(L^p)^*=\{0\}$ when $p<1$ and
$\mathcal{M}$ has no minimal projection; compare that with

\begin{corollary}
If $\mathcal{M}$ has no minimal projection and $p<q$,
$$\text{Hom}(L_\mathcal{M}^p,L_\mathcal{M}^q) = \{0\}.$$
\end{corollary}
\begin{proof}
Choose a state $\varphi$. ÊIf $T$ is a bounded morphism and
$\frac{1}{p} =\frac{1}{q} + \frac{1}{r}$, set
$$\widetilde{T}:L^q_\mathcal{M} \to L^q_\mathcal{M} \text{ by }
\widetilde{T}(\xi) \mapsto T(\varphi^{1/r} \xi).$$
This is a bounded module map, so by the preceding corollary there must be $x \in
\mathcal{M}$ with
$$x\xi = \widetilde{T}(\xi) = T(\varphi^{1/r} \xi).$$
If $x \ne 0$, let $x = v|x|$ and $e = e(\varepsilon, \infty)$ be a nonzero spectral projection of $|x|$. ÊFor all $\xi \in L^p$, we have
$$\|e\xi\|_q = \|ve\xi\| = \|v |x| e |x|^{-1} e \xi\| = \|T(\varphi^{1/r} e
|x|^{-1} e \xi) \| =\| T(\eta e \xi) \| \le C \|\eta e \xi \|_p,$$
where $\eta = \varphi^{1/r} e |x|^{-1} e.$ ÊIt remains to show that such a
``reversed H\"{o}lder inequality" cannot hold.

Let $f_n$ be a decreasing sequence of nonzero projections $\le e$ and converging
strongly to 0. Ê(This is where nonatomicity is essential.) ÊThen by Lemma
$\ref{T:strong}$, $\| \eta f_n \| \to 0.$ ÊChoose an element $f$ with
$$\|\eta f \|_r < 1/C.$$
Now take a functional $\rho$ with $s(\rho) = f.$ ÊIt follows that
$$\|\rho^{1/q}\|_q \le C \| \eta f \rho^{1/q} \|_p \le C \|\eta f\|_r \|
\rho^{1/q} \|_q < \|\rho^{1/q}\|_q,$$
which is impossible. ÊSo $x = 0$, which implies $T(\varphi^{1/r} \cdot)$ is the zero map. ÊSince this holds for any choice of $\varphi$, $T$ must also be the zero map.
\end{proof}

\section{$L^p$ modules}

Now we turn to the development of an $L^p$ representation theory. ÊNote that
this cannot mean representations on classical $L^p$ spaces: $L^p(\mathcal{M})$
itself is not a classical $L^p$ space unless $p=2$ or $\mathcal{M}$ is
commutative. ÊWe would like to build the category out of $L^p(\mathcal{M})$ in
the same way that nondegenerate normal right Hilbert space representations are
built out of $L^2(\mathcal{M})$.

Let us examine a countably generated Hilbert module $\mathfrak{H}_\mathcal{M}$. ÊFollowing standard arguments, $\mathfrak{H}$ decomposes into a direct sum of cyclic representations $(\overline{\xi_n \mathcal{M}})_\mathcal{M}$, each of which is isomorphic to the GNS representation for the associated vector state, and all GNS representations are reductions of $L^2(\mathcal{M})$. ÊSo we have
$$\mathfrak{H}_\mathcal{M} \simeq (\bigoplus \overline{\xi_n
\mathcal{M}})_\mathcal{M} \simeq (\bigoplus
\mathfrak{H}_{\omega_{\xi_n}})_\mathcal{M} \simeq (\bigoplus q_n
L^2(\mathcal{M}))_\mathcal{M}.$$
(In fact $q_n =s_\ell(\xi_n).$)

Since this is a right module, it is natural to write vectors as columns with the $n$th entry in $q_n L^2$:
$$\mathfrak{H} \simeq \left(
\begin{smallmatrix}
q_1 L^2(\mathcal{M}) \\
q_2 L^2(\mathcal{M}) \\
\vdots
\end{smallmatrix}
\right) \simeq (\sum q_n \otimes e_{nn}) \left(
\begin{smallmatrix}
L^2(\mathcal{M}) \\
L^2(\mathcal{M}) \\
\vdots
\end{smallmatrix}
\right).$$
Here $e_{nn}$ are diagonal matrix units in $M_\infty$, so $(\sum q_n \otimes e_{nn})$ is a diagonal projection in $\mathcal{M} \otimes \mathcal{B}(\mathfrak{H}).$ ÊThe right action of $\mathcal{M}$ is, of course, matrix multiplication (by $1\times 1$ matrices) on the right. ÊModules which are not countably generated can be represented by columns and projections of larger size, and non-diagonal projections work equally well - see Section 5.

Our target class of modules is obtained by replacing the index 2 by $p$.
Although this seems simple enough, the geometry of such spaces presents certain
difficulties. ÊTo start with, one cannot obtain the norm of a column via an
$\ell^p$ (or $\ell^2$) sum. ÊThe following example will serve as motivation.

Consider the right $L^p$ $\mathcal{M}$-module
$$ \mathfrak{X} =\left(
\begin{smallmatrix}
L^p \\
L^p
\end{smallmatrix}
\right).$$
This should be a left $L^p$ $M_2(\mathcal{M})$-module, as it is
$$L^p(M_2(\mathcal{M})) e_{11} =
\left(
\begin{smallmatrix}
L^p(\mathcal{M}) & 0 \\
L^p(\mathcal{M}) & 0
\end{smallmatrix}
\right).$$
It is then a left submodule of $L^p(M_2(\mathcal{M}))$ and so inherits the norm:
$$\Vert
\left(
\begin{smallmatrix}
\xi \\
\eta
\end{smallmatrix}
\right)
\Vert =
\Vert
\left(
\begin{smallmatrix}
\xi Ê& 0 \\
\eta & 0
\end{smallmatrix}
\right)
\Vert_{L^p(M_2(\mathcal{M}))} = \|\xi^* \xi + \eta^* \eta\|_{p/2}^{1/2},$$
which is in general not purely a function of the norms of $\xi$ and $\eta$.

Norm-determining expressions of the form $\xi^* \xi$ recall inner products in
Hilbert C*-modules. ÊBased on this parallel, we make

\begin{definition}

Let $\mathfrak{X}$ be a complex vector space which is a right
$\mathcal{M}$-module and $p \in (0, \infty]$. ÊBy an \textbf{$L^{p/2}$-valued
inner product} on $\mathfrak{X}$ we mean a sesquilinear mapping, conjugate
linear in the first variable, from $\mathfrak{X} \times \mathfrak{X}$ to
$L^{p/2}(\mathcal{M})$ which satisfies
\begin{enumerate}
\item[(i)] $<\xi, \eta x> = <\xi, \eta> x;$
\item[(ii)] $<\xi, \eta> = <\eta, \xi>^*;$
\item[(iii)] $<\xi, \xi> \ge 0; \quad <\xi, \xi> = 0 \iff \xi = 0.$
\end{enumerate}
\end{definition}

\begin{proposition}
\begin{equation}\label{E:middle}
<\xi, \eta> = <\xi, \xi>^{1/2} T <\eta, \eta>^{1/2}
\end{equation}
for some $T \in \mathcal{M}$ with $\|T\| \le 1.$ ÊSo if we set
$$||\xi|| \triangleq ||<\xi, \xi>^{1/2}||_p,$$
then
\begin{equation} \label{E:cs}
\|<\xi, \eta>\|_{p/2} \le \|\xi\| \|\eta\|.
\end{equation}
We have that $\| \cdot \|$ is a norm when $p \ge 2$ and a $p/2$-norm when $p \le
2$. Ê(This is improved by the end of the next section.)
\end{proposition}
\begin{proof}
Most of this proof is standard. ÊFor $\xi, \eta \in \mathfrak{X}$, consider the
matrix
$$A = \left(
\begin{smallmatrix}
<\xi, \xi> & <\xi, \eta>\\
<\eta, \xi> & <\eta, \eta>
\end{smallmatrix}
\right) \in M_2(L^{p/2}(\mathcal{M})) \simeq L^{p/2}(M_2(\mathcal{M})).$$
We claim that $A$ is positive. ÊIf $\mathcal{M}$ is semifinite, we may choose a faithful semifinite trace $\tau$ and consider the $L^p$ spaces to be spaces of $\tau$-measurable operators. ÊFor $x,y \in L^2 \cap L^\infty$, which is dense in $L^2$,
$$\left< \left( \begin{smallmatrix}
<\xi, \xi> & <\xi, \eta>\\
<\eta, \xi> & <\eta, \eta>
\end{smallmatrix} \right)
\left( \begin{smallmatrix} x \\ y \end{smallmatrix} \right), Ê\left(
\begin{smallmatrix} x \\ y \end{smallmatrix} \right) \right> = <\xi x + \eta y,
\xi x + \eta y> \;\ge 0,$$
and by density the matrix is positive.

If $\mathcal{M}$ is purely infinite, then so is $M_2(\mathcal{M})$; let $v$ be a
partial isometry in $M_2(\mathcal{M})$ with $v v^* = 1$, $v^*v = e_{11}$. ÊThus
$v$ is of the form
$$\left( \begin{smallmatrix} v_{11} & 0 \\ v_{12} & 0 \end{smallmatrix}
\right).$$
We have
\begin{align*}
A &= v v^* A v v^* = v \left( \begin{smallmatrix} v_{11}^* & v_{12}^* \\ 0 & 0
\end{smallmatrix} \right) A \left( \begin{smallmatrix} v_{11} & 0 \\ v_{12} & 0
\end{smallmatrix} \right) v^* \\ &= v \left( \begin{smallmatrix} <\xi v_{11} +
\eta v_{12}, \xi v_{11} + \eta v_{12}> & 0 \\ 0 & 0 \end{smallmatrix} \right)
v^* \ge 0.
\end{align*}
A von Neumann algebra decomposes as a direct sum of semifinite and purely
infinite summands, so we see that $A$ is positive in general.

Now the usual matrix manipulations give \eqref{E:middle}, and \eqref{E:cs}
follows by H\"{o}lder's inequality.

When $p \le 2$, use the inequality from [K2]
$$\|v+w\|_q^q \le \|v\|_q^q + \|w\|_q^q, \quad v,w \in L^q, \quad q\le 1$$
to write
\begin{align} \label{E:pnorm}
\|\xi + \eta\|^p &= \|<\xi, \xi> + <\xi, \eta> + <\eta, \xi> + <\eta,
\eta>\|_{p/2}^{p/2} \\
\notag &\le \|<\xi, \xi>\|_{p/2}^{p/2} Ê+ \|<\xi, \eta>\|_{p/2}^{p/2} Ê+ \|<\eta,
\xi>\|_{p/2}^{p/2} Ê+ \|<\eta, \eta>\|_{p/2}^{p/2} \\
\notag &\le \|\xi\|^p +
2\|\xi\|^{p/2} \|\eta\|^{p/2} + \|\eta\|^p \\ \notag &= (\|\xi\|^{p/2} +
\|\eta\|^{p/2})^2.
\end{align}
Therefore
$$\| \xi + \eta \|^{p/2} \le \| \xi\|^{p/2} + \| \eta \|^{p/2}.$$

When $p \ge 2$, one starts $\eqref{E:pnorm}$ with $\| \xi + \eta \|^2$ and
proves the triangle inequality via the same manipulations.
\end{proof}

It is worth noting that $\|\xi x \| \le \|\xi\| \|x\|,$ so the action of
$\mathcal{M}$ is continuous.

\begin{definition} \label{D:lpmod}
For $p < \infty$, a right $\mathcal{M}$-module $\mathfrak{X}$
is called a \textbf{right $L^p$ $\mathcal{M}$-module} if it has an
$L^{p/2}$-valued inner product and is complete in the inherited (quasi)norm.
For $p=\infty$, we keep this condition (so $\mathfrak{X}$ is a Hilbert $C^*$
module) and impose the additional requirement that the unit ball of
$\mathfrak{X}$ be closed in the strong topology, i.e. the topology arising from
the seminorms
$$\xi \mapsto ( \varphi(<\xi,\xi>) )^{1/2}, \quad \varphi \in \mathcal{M}_*^+.$$
\end{definition}

The set $\overline{<\mathfrak{X}, \mathfrak{X}>}$ is a closed self-adjoint sub-bimodule of $L^{p/2}$, which must have the form $zL^{p/2}$ for some central projection $z \in \mathcal{M}$. ÊSo $\mathfrak{X}$ is a \textit{faithful} right $L^p$ $z\mathcal{M}$-module.

\textbf{Examples:}
\begin{itemize}
\item Any classical $L^p$ space is a right $L^p$ module for the corresponding
$L^\infty$ algebra, with inner product $$<f,g> = \bar{f}g,$$.

\item Any normal right representation of $\mathcal{M}$ on a Hilbert space
$\mathfrak{H}=\mathfrak{X}$ admits a unique structure as $L^2$-module by setting
the inner product $$<\xi,\eta>_{\mathfrak{X}}$$ to be the state $\omega_{\xi,
\eta}$ defined by $$<\xi,\eta x>_{\mathfrak{H}} = \omega_{\xi, \eta}(x).$$ Ê(For
coherence, the inner product in $\mathfrak{H}$ should be linear in the second
argument.) ÊOn the other hand, any $L^2$ $\mathcal{M}$-module $\mathfrak{X}$ is
also a Hilbert space via $$<\xi, \eta>_{\mathfrak{H}} \triangleq
<\xi,\eta>_{\mathfrak{X}} (1) = \text{Tr}(<\xi,\eta>_{\mathfrak{X}}),$$ where Tr
denotes the Haagerup trace on $L^1$. ÊSince the $\mathfrak{X}$-inner product is $\mathcal{M}_*$-valued, the Hilbert space representation is automatically normal.

\item $L^p(\mathcal{M})$ is a right $L^p$ module with inner product $<\xi,\eta>
= \xi^* \eta.$ ÊSimilarly for $qL^p$, where $q$ is a projection in
$\mathcal{M}$.

\end{itemize}

We wish to highlight a special class of right $L^p$ $\mathcal{M}$-modules; call
them `principal' for the time being. ÊIf $\{q_\alpha\}_{\alpha \in I}$ are
projections in $\mathcal{M}$, the set
\begin{equation}
\left\{(\xi_\alpha)\;|\;\xi_\alpha \in q_\alpha L^p,\; \sum \xi_\alpha^*
\xi_\alpha \in L^{p/2}\right\}
\end{equation}
Êis a right $L^p$ $\mathcal{M}$-module with $$<(\xi_\alpha), (\eta_\alpha)> =
\sum \xi_\alpha^* \eta_\alpha.$$
For $p=\infty$, Paschke [P] showed that the directed net of finite sums
converges strongly; he called this construction an \textit{ultraweak direct
sum}.

For $p<\infty$, the limit (of finite sums) exists in norm. ÊThis follows from
the Cauchy-Schwarz inequality, which can be proven directly as follows. ÊLet $\mathfrak{H}_I$
be the Hilbert space with dimension $|I|$, and set
$$\tilde{\xi} = \sum \xi_\alpha \otimes e_{\alpha 1} \in L^p(\mathcal{M} \otimes
\mathcal{B}(\mathfrak{H}_I)).$$
(So we are placing $\xi$ along the first column of a matrix.) ÊThen Kosaki's
generalized H\"{o}lder inequality [K] guarantees that
$$\tilde{\xi}^*\tilde{\eta} \in L^{p/2}(\mathcal{M} \otimes
\mathcal{B}(\mathfrak{H}_I))$$
and
\begin{align*}
||<(\xi_\alpha), (\eta_\alpha)>||_{p/2} &= ||\tilde{\xi}^*\tilde{\eta}||_{p/2} \le
||\tilde{\xi}||_p ||\tilde{\eta}||_p = \left\|(\sum \xi_\alpha^*
\xi_\alpha)^{1/2}\right\|_p \left\|(\sum \eta_\alpha^*
\eta_\alpha)^{1/2}\right\|_p \\ &=
\left\|(\sum \xi_\alpha^* \xi_\alpha)\right\|_{p/2}^{1/2} \left\|(\sum
\eta_\alpha^* \eta_\alpha)\right\|_{p/2}^{1/2} = ||\xi|| ||\eta||.
\end{align*}

We denote this module $\bigoplus_c q_\alpha L^p$ for \textit{column} sum.
Indeed, the reader should think of principal modules as columns with
entries from $L^p$. ÊMotivated by this, we make

\begin{definition}
Let $\{\mathfrak{X}_\alpha \}$ be $L^p$ modules. ÊIf $(\xi_\alpha)$ and
$(\eta_\alpha)$ have finite support, set
$$<(\xi_\alpha), (\eta_\alpha)> = \sum < \xi_\alpha, \eta_\alpha >.$$
The \textbf{column sum}, $\bigoplus_c \mathfrak{X}_\alpha$, is the closure of
the finitely supported vectors with respect to the (quasi)norm ($p < \infty$) or
strong topology ($p = \infty$) coming from this inner product.

We denote the countable column sum of $L^p(\mathcal{M})$ as $C^p(\mathcal{M})$, or simply $C^p$ if the underlying algebra is clear.

Note: As above, it will turn out that $\bigoplus_c \mathfrak{X}_\alpha = \{
(\xi_\alpha ) | \sum_\alpha <\xi_\alpha, \xi_\alpha> \in L^{p/2} \}$.

\end{definition}

We can now state one of our main results.
\begin{theorem}\label{T:main}
Any $L^p$ module is isometrically isomorphic, as a module, to a principal $L^p$
module.
\end{theorem}

If $\mathfrak{X}$ is cyclic, this is easy. ÊTake
$$\mathfrak{X} = \overline{\xi \mathcal{M}}$$
and consider the densely-defined isomorphism of $L^p$ modules
$$\mathfrak{X} \leftrightarrow s(<\xi, \xi>^{1/2}) L^p: \quad \xi x
\leftrightarrow <\xi, \xi>^{1/2} x.$$
Since the inner product and the bounded action of $\mathcal{M}$ extend
continuously to the completion, this is an isomorphism.

The whole difficulty of the proof lies in devising the column sum decomposition. ÊThis may be thought of as a generalization of the fact that Hilbert spaces have an orthonormal basis. Ê(A version of this theorem was proven for a special type of $L^p$ module in [J1, Prop. 2.8].)

\section{Proof of Theorem \ref{T:main}}
If $p=2$, $\mathfrak{X}$ is a Hilbert space. ÊThe previously mentioned decomposition theorem gives
$$\mathfrak{X} \simeq \bigoplus q_\alpha L^2(\mathcal{M})$$
isometrically as modules. ÊNow the right-hand side admits a unique $L^1$-valued
inner product and so is a column sum in our sense; therefore $\mathfrak{X}$ is
principal.

We consider the cases $p=\infty$, $p>2$, and $p<2$ separately.

\underline{\textit{Case 1: $p=\infty$}}

Choose $\xi \in \mathfrak{X}$ and set
$$\xi_n = \xi \left( <\xi, \xi> + \frac{1}{n} \right)^{-1/2} \in \mathfrak{X},$$
so $\| \xi_n\| < 1.$

We observe
$$\varphi(<\xi_m - \xi_n, \xi_m - \xi_n>) \to 0, \quad m,n \to 0.$$
Then $\xi_n$ converges strongly, say to $\eta_1$, and apparently $<\eta_1,
\eta_1>$ is a projection $q_1$.

Consider a maximal set $\{ \eta_\alpha \}$ with the property that
$$<\eta_\beta, \eta_\gamma> = \delta_{\beta \gamma} q_\beta.$$
If the strong closure of $\sum \eta_\alpha \mathcal{M}$ is not all of
$\mathfrak{X}$, choose $\xi$ outside this set and write
\begin{equation} \label{E:bessel}
\xi = \sum \eta_\alpha <\eta_\alpha, \xi> + \left(\xi - \sum \eta_\alpha
<\eta_\alpha, \xi> \right).
\end{equation}
The first summand should be interpreted as a strong limit; existence follows
from the Bessel-type inequality
$$0 \le \left< \xi - \sum \eta_\alpha <\eta_\alpha, \xi>, \: \xi - \sum \eta_\alpha <\eta_\alpha, \xi> \right> = <\xi, \xi> - \sum |<\eta_\alpha, \xi>|^2.$$
By assumption the second summand in \eqref{E:bessel} is nonzero. ÊWe
can normalize it as above (which does not change orthogonality) and add it to
our set $\{ \eta_\alpha \}$ - but this violates the maximality of $\{
\eta_\alpha \}$.

Therefore the strong closure of $\sum \eta_\alpha \mathcal{M}$ is
$\mathfrak{X}$. ÊFinally we have an isomorphism
$$\mathfrak{X} \ni \xi \leftrightarrow (<\eta_\alpha, \xi>) \in \bigoplus_c
q_\alpha \mathcal{M}.$$

Essentially this is Paschke's result [P], but we have started with a topological
condition instead of an algebraic one (self-duality). ÊA Hilbert C*-module
$\mathfrak{X}$ is called \textit{self-dual} if $\mathfrak{X} \simeq \text{Hom} (
\mathfrak{X}_\mathcal{M}, \mathcal{M}_\mathcal{M} )$ via $\xi \leftrightarrow <
\xi, \cdot \,>.$ ÊWeaker than the strong topology we have defined is the
\textit{weak} topology on the unit ball, generated by the functionals
$$\xi \mapsto \varphi(<\eta, \xi>), \quad \varphi \in \mathcal{M}_*, \: \eta \in
\mathfrak{X}.$$

We have arrived at
\begin{theorem}
For a Hilbert C*-module $\mathfrak{X}$ over a von Neumann algebra $\mathcal{M}$,
the following conditions are equivalent:
\begin{itemize}
\item[(i)] the unit ball of $\mathfrak{X}$ is strongly closed;
\item[(ii)] $\mathfrak{X}$ is principal; or, to say the same thing,
$\mathfrak{X}$ is an ultraweak direct sum of Hilbert C*-modules $q_\alpha
\mathcal{M}$, for some projections $q_\alpha$;
\item[(iii)] $\mathfrak{X}$ is self-dual;
\item[(iv)] the unit ball of $\mathfrak{X}$ is weakly closed.
\end{itemize}
\end{theorem}

This theorem has consequences for an arbitrary Hilbert C*-module $\mathfrak{X}$
over a von Neumann algebra $\mathcal{M}$. ÊSet $\bar{\mathfrak{X}}$ to be the
strong closure; a straightforward argument shows that $\bar{\mathfrak{X}}$ is an
$L^\infty$ module for $\mathcal{M}$. ÊTherefore $\mathfrak{X}$ is representable
as a strongly dense submodule of a principal $L^\infty$ module. ÊThis
observation, and a similar discussion, are also found in [We].

\underline{\textit{Case 2: $p>2$}}

Let $\{\xi_\alpha\}$ be a maximal orthogonal set (with no condition on
$<\xi_\alpha, \xi_\alpha>$). ÊSet
$$\mathfrak{X}_0 = \overline{\sum \overline{\xi \mathcal{M}}} \simeq \bigoplus_c q_\alpha L^p.$$
Any vector $\eta$ in $\mathfrak{X}_0$ can be written as a limit, i.e.
$$ \eta = \lim_n \sum \xi_\alpha x_{\alpha, n}.$$
But if this is Cauchy, the orthogonality of $\{\xi_\alpha\}$ implies
\begin{align*}
0 &= \lim_{m,n} \|\sum \xi_\alpha x_{\alpha, n} - \sum \xi_\alpha x_{\alpha, m}\| \\ &\ge \lim_{m,n} \| \xi_\alpha (x_{\alpha, n} - x_{\alpha, m}) \| \quad
\text{ for each $\alpha$.}
\end{align*}
Thus $\eta$ has a unique representation as $\sum \eta_\alpha$, $\eta_\alpha \in
q_\alpha L^p,$ and we have an isomorphism of $L^p$ modules
$$\mathfrak{X}_0 \ni \eta \leftrightarrow \bigoplus_c q_\alpha L^p.$$

So we just need to show that $\mathfrak{X}_0 = \mathfrak{X}.$ ÊNow
$\mathfrak{X}$ is a Banach space since $p>2$, and $\mathfrak{X}_0$ was seen to
be reflexive in the last paragraph. ÊTherefore $\mathfrak{X}_0$ is a proximinal
subspace of $\mathfrak{X}$ ([Si], Cor. 2.1). ÊThis means that if $\xi \in
\mathfrak{X} \setminus \mathfrak{X}_0$, there exists an element $\eta_0$ in
$\mathfrak{X}_0$ with
$$\|\xi - \eta_0 \| = \inf_{\eta \in \mathfrak{X}_0} \|\xi - \eta \|.$$
Then $\zeta = \xi - \eta_0$ has 0 as a best approximant. ÊWhat does this say
about $\zeta$?

Fix $\alpha$ and $x \in \mathcal{M}$. ÊBy assumption, the function
$$\mathbb{R} \ni t \mapsto \| \zeta + t \xi_\lambda x \|$$
attains its minimum at $t=0$. ÊSet
$$<\zeta, \xi_\lambda> = <\zeta, \zeta>^{1/2} T_\lambda <\xi_\lambda,
\xi_\lambda>^{1/2},$$
using $\eqref{E:middle}$, and observe
\begin{align*}
\| \zeta &+ t \xi_\lambda x \| \\ &= \|<\zeta, \zeta> + 2t \text{Re }(<\zeta,
\zeta>^{1/2} T_\lambda <\xi_\lambda, \xi_\lambda>^{1/2} x) + t^2
x^*<\xi_\lambda, \xi_\lambda> x\| \\ &\overset{o(t)}{\thicksim} \|<\zeta,
\zeta> + 2t \text{Re }(<\zeta, \zeta>^{1/2} T_\lambda <\xi_\lambda, \xi_\lambda
>^{1/2} x) \\ &{} \qquad + t^2 x^*<\xi_\lambda, \xi_\lambda>^{1/2} T_\lambda^* T_\lambda
<\xi_\lambda, \xi_\lambda>^{1/2}x \| \\ &= \|<\zeta, \zeta>^{1/2} + t T_\lambda
<\xi_\lambda, \xi_\lambda>^{1/2} x \| \triangleq f_{\lambda, x}(t).
\end{align*}
Now $ f_{\lambda, x}$ is differentiable since $\zeta$ was presumed nonzero and
the norm in $L^p \setminus\{0\}$ is Fr\'{e}chet differentiable. Ê(The Clarkson inequalities imply $L^p$ is uniformly smooth for
$1<p<\infty$ [QX].) ÊIt agrees up to $o(t)$ with a function which has a local minimum
at $t=0$, so $ f_{\lambda, x}'(0)=0$. ÊFinally, it is convex by construction.
It follows that $ f_{\lambda, x}$ attains its absolute minimum ($=\|<\zeta,
\zeta>^{1/2}\|$ $ = \| \zeta \|)$ at 0.

Since this is true for all $\lambda$ and $x$ we get that in $L^p$
$$\text{dist}(<\zeta, \zeta>^{1/2}, \overline{\sum T_\lambda <\xi_\lambda,
\xi_\lambda>^{1/2} \mathcal{M}}) = \|\zeta\|.$$
By Hahn-Banach there is a norm one functional on $L^p$ which annihilates the
subspace $\overline{\sum T_\lambda <\xi_\lambda, \xi_\lambda>^{1/2}
\mathcal{M}}$ and takes the value $\| \zeta \|$ at $<\zeta,\zeta>^{1/2}$. ÊThis
functional must have the form $\text{Tr}(v\rho^{1/q} \cdot)$ for some $v \in
\mathcal{M}$, $\rho \in \mathcal{M}_*^+.$ ÊThen we fix $\lambda$ and write out
$$\text{Tr}(v \rho^{1/q} T_\lambda <\xi_\lambda, \xi_\lambda>^{1/2} \mathcal{M})
= 0$$
\begin{equation}\label{E:rho1}
\Rightarrow s(\rho) \bot s_\ell(T_\lambda <\xi_\lambda, \xi_\lambda>^{1/2}).
\end{equation}
Also
$$\|<\zeta, \zeta>^{1/2}\|_p = |\text{Tr}(v \rho^{1/q} <\zeta, \zeta>^{1/2})|
\le \|s(\rho) <\zeta, \zeta>^{1/2}\|_p \le \|<\zeta, \zeta>^{1/2} \|_p,$$
so these are equalities and in particular
\begin{equation}\label{E:rho2}
s(\rho) = s(<\zeta, \zeta>^{1/2}).
\end{equation}
Together \eqref{E:rho1} and \eqref{E:rho2} imply
$$<\zeta, \xi_\lambda> = <\zeta, \zeta>^{1/2} T_\lambda <\xi_\lambda,
\xi_\lambda>^{1/2} = 0.$$
So the set $\{\xi_\lambda\}$ was not a maximal orthogonal set in $\mathfrak{X}$,
a contradiction. ÊThis completes the proof for $p>2.$

\underline{\textit{Case 3: $p<2$}}

By restricting the algebra (see the discussion following Definition \ref{D:lpmod}), we may assume that the module is faithful. ÊWe need two auxiliary constructions.

\textbf{I.} Let $\mathfrak{X}$ be an $L^p$ module, and $1/p + 1/q = 1/r$. ÊWe
write
$$\mathfrak{X} \otimes_{\mathcal{M}} L^q(\mathcal{M})$$
for the closure of the algebraic tensor product, modulo the null space, in the topology arising from the degenerate inner product
$$<\xi_1 \otimes \eta_1, \xi_2 \otimes \eta_2> = \eta_1^* <\xi_1, \xi_2> \eta_2
\in L^{r/2}(\mathcal{M}).$$
It is easy to see that $\mathfrak{X} \otimes_{\mathcal{M}} L^q(\mathcal{M})$
satisfies the relation
\begin{equation}
\xi x \otimes \eta = \xi \otimes x \eta
\end{equation}
and is an $L^r$ module in our sense.

\textbf{II.} Let $\mathfrak{X}$ be an $L^p$ module and $\varphi$ be a fixed
faithful \textit{strictly semifinite} weight on $\mathcal{M}$. ÊThis means that $\varphi = \sum \varphi_\alpha$, where the $\varphi_\alpha$ are orthogonal and bounded. ÊWe will create an $L^2$ module with the same ``shape" as
$\mathfrak{X}$.

\begin{lemma}\label{L:Dphi}
The following conditions on a vector $\xi \in \mathfrak{X}$ are equivalent:
\begin{enumerate}
\item[(i)] $<\xi, \xi> \le C\varphi^{2/p}$ for some $C$;
\item[(ii)] $<\xi, \xi>^{1/2} = y \varphi^{1/p}$ for some $y \in \mathcal{M}$;
\item[(iii)] $<\xi, \xi>^{1/2} = \varphi^{1/p} z$ for some $z \in \mathcal{M}$.
\end{enumerate}
We denote the set of such vectors as $\mathcal{D}_\varphi$.
\end{lemma}
This is nothing but \eqref{E:order}.

\begin{lemma} \label{T:phidense}
$\mathcal{D}_\varphi$ is dense in $\mathfrak{X}.$
\end{lemma}
\begin{proof}
Given any $\xi \in \mathfrak{X}$, let $q = s(<\xi, \xi>)$ and densely define the $L^p$ module isomorphism $T$ by
$$T: q L^p \overset{\sim}{\to} \overline{\xi \mathcal{M}} \subset \mathfrak{X}, \qquad T(<\xi, \xi>^{1/2} x) = \xi x.$$

We need that elements of the form $y \varphi^{1/p}$ are dense in $L^p$. ÊBecause $\varphi^{1/p}$ is not necessarily $\tau$-measurable, this is slightly more delicate than Lemma \ref{T:dense}.

Set $q_\alpha = s(\varphi_\alpha)$, and let $\{r_\beta\}$ be the net of finite sums of the $q_\alpha$ (ordered naturally). ÊAgain by Lemma \ref{T:addcols}, the net $\{<\xi, \xi>^{1/2} r_\beta \}$ converges to $<\xi, \xi>^{1/2}.$ ÊSince $r_\beta$ commutes with $\varphi^{1/p}$, we have
$$<\xi, \xi>^{1/2} r_\beta \in L^p r_\beta = \overline{\mathcal{M} r_\beta \varphi^{1/p}}.$$
Putting these two approximations together, we may find $\{y_n\} \subset \mathcal{M}$ with $y_n \varphi^{1/p} \to$ $ <\xi, \xi>^{1/2}$.

$T$ is an isomorphism and preserves inner products, so
$$T(y_n \varphi^{1/p}) \to T(<\xi, \xi>^{1/2}) = \xi.$$
Also
$$< T(y_n \varphi^{1/p}), T(y_n \varphi^{1/p})> = < y_n \varphi^{1/p}, y_n \varphi^{1/p}> = \varphi^{1/p}y_n^* y_n \varphi^{1/p} \le \|y\|^2 \varphi^{2/p},$$
so by Lemma \ref{L:Dphi}, $ T(y_n \varphi^{1/p}) \in \mathcal{D}_\varphi$. ÊThus we have written $\xi$ as a limit of vectors in $\mathcal{D}_\varphi$.

\end{proof}

With $\frac{1}{2} + \frac{1}{r} = \frac{1}{p}$, we define an $L^1$-valued inner
product on $\mathcal{D}_\varphi$ by
\begin{equation}\label{E:Dinn}
<\xi, \eta>_{\mathcal{D}_\varphi} \triangleq
\varphi^{-1/r}<\xi,\eta>_\mathfrak{X} \varphi^{-1/r}.
\end{equation}
By \eqref{E:middle} and ÊLemma \ref{L:Dphi}, $\varphi^{1/p}$ factors out of $<\xi, \eta>$ on both the left and the right, and \eqref{E:Dinn} is justified. ÊThe nontrivial fact that composition with $\varphi^{-1/r}$ is the inverse of composition with $\varphi^{1/r}$ is found in [S2].

We now describe the module action. ÊClearly the previous $\mathcal{M}$-action is
not compatible with the new inner product (and $\mathcal{D}_\varphi$ is not a
submodule of $\mathfrak{X}$). ÊInstead we need to work with
$\mathcal{M}^\varphi_a$, the operators in $\mathcal{M}$ for which
$$t \mapsto \sigma^\varphi_t(x) = \varphi^{it} x \varphi^{-it}$$
extends off the real line to an entire $\mathcal{M}$-valued function.

The action must be
$$\eta \cdot x = \eta \varphi^{-1/r} x \varphi^{1/r};$$
then
\begin{align*}
<\xi, \eta \cdot x>_{\mathcal{D}_\varphi} &= <\xi, \eta \varphi^{-1/r} x
\varphi^{1/r}>_{\mathcal{D}_\varphi} = \varphi^{-1/r}<\xi, \eta \varphi^{-1/r} x
\varphi^{1/r}>_{\mathfrak{X}} \varphi^{-1/r} \\ &= \varphi^{-1/r}<\xi,
\eta>_{\mathfrak{X}} \varphi^{-1/r} x \varphi^{1/r} \varphi^{-1/r} = <\xi,
\eta>_{\mathcal{D}_\varphi} x.
\end{align*}

As we noted before, an $L^1$-valued inner product composed with Tr is a usual
inner product; therefore the closure of $\mathcal{D}_\varphi$ in the inner
product norm is a Hilbert space $\mathfrak{H}_{\mathfrak{X},\varphi}$. ÊThe
*-algebra $\mathcal{M}^\varphi_a$ is represented isometrically on it - in fact
it is a *-representation:
\begin{align*}
<\xi, \eta \cdot x^*>_{\mathfrak{H}_{\mathfrak{X},\varphi}} &= \text{Tr} (<\xi,
\eta \cdot x^*>_{\mathcal{D}_\varphi}) = \text{Tr} (<\xi,
\eta>_{\mathcal{D}_\varphi} x^*) \\ &= \text{Tr} (x^*<\xi,
\eta>_{\mathcal{D}_\varphi}) = \text{Tr} (<\xi \cdot x,
\eta>_{\mathcal{D}_\varphi}).
\end{align*}

We need to show that the von Neumann closure of $\mathcal{M}^\varphi_a$ is
exactly $\mathcal{M}$. ÊA dense set of vector states in this representation is
\begin{equation} \label{E:vectors}
x \mapsto <\xi, \xi \cdot x>_{\mathfrak{H}_{\mathfrak{X},\varphi}}, \quad \xi
\in \mathcal{D}_\varphi,
\end{equation}
and these are identical to the linear functionals
$$x \mapsto \text{Tr} (<\xi, \xi>_{\mathcal{D}_\varphi} x), \quad \xi \in
\mathcal{D}_\varphi.$$
Deducing further and using Lemma \ref{L:Dphi},
\begin{align*}
\{<\xi, \xi>_{\mathcal{D}_\varphi} \mid \xi \in \mathcal{D}_\varphi \} &= \{\varphi^{-1/r} \psi^{2/p} \varphi^{-1/r} \mid C\varphi^{2/p} \ge \psi^{2/p} \in L^{p/2} \}\\ &= \{\varphi^{1/2} |y|^2 \varphi^{1/2} \mid y\varphi^{1/p} \in L^p \}.
\end{align*}
Now we need another double approximation argument, and we are brief. ÊSince $\varphi$ is semifinite, any element of $\mathcal{M}_*^+$ is a norm limit of elements $\varphi^{1/2} |y_n|^2 \varphi^{1/2}$, where $y_n \in \mathfrak{N}_\varphi$, the definition ideal of $\varphi$. ÊEach of these can be approximated by an element
$$r_\beta \varphi^{1/2} |y|^2 \varphi^{1/2} r_\beta = \varphi^{1/2} r_\beta |y|^2 r_\beta \varphi^{1/2}$$
($r_\beta$ are as in the proof of Lemma \ref{T:phidense}), and these belong to the sets above since $y r_\beta \varphi^{1/p} \in L^p.$

The upshot of all this is that the vector states in \eqref{E:vectors} form a dense set in $\mathcal{M}_*^+.$ ÊThus the strong topology in this representation agrees with the strong topology in the representation of
$\mathcal{M}^\varphi_a$ on $L^2(\mathcal{M})$. ÊHappily, $\mathcal{M}^\varphi_a$
is dense in $\mathcal{M}$ in the latter topology, so the von Neumann closure is
$\mathcal{M}$.

The reader can check that the extensions of the $\mathcal{M}$-action and
$L^1$-valued inner product to $\mathfrak{H}_{\mathfrak{X},\varphi}$ do make it
into an $L^2$ module for $\mathcal{M}$.

\bigskip

Now consider the $L^p$ module
$$\mathfrak{H}_{\mathfrak{X},\varphi} \otimes_\mathcal{M} L^r.$$
We will make two observations: that it is principal, and that it is isomorphic
to $\mathfrak{X}$.

$\mathfrak{H}_{\mathfrak{X},\varphi}$ is an $L^2$ module and so of the form
$\bigoplus_c q_\alpha L^2$. ÊIt is not hard to see that the functor
``$\otimes_\mathcal{M} L^r$" commutes with column sums; i.e.
$$\mathfrak{H}_{\mathfrak{X},\varphi} \otimes_\mathcal{M} L^r = \left(
\bigoplus_c q_\alpha L^2 \right) \otimes_\mathcal{M} L^r \simeq \bigoplus_c
\left( q_\alpha L^2 \otimes_\mathcal{M} L^r \right) = \bigoplus_c (q_\alpha
L^p),$$
which is principal.

Consider the dense submodule
$$\mathcal{D}_\varphi \otimes_\mathcal{M} (\varphi^{1/r} \mathcal{M} \cap L^r) \subset
\mathfrak{H}_{\mathfrak{X},\varphi} \otimes_\mathcal{M} L^r.$$
For elements of this subset, we have
$$<\xi \otimes \varphi^{1/r}x, \eta \otimes \varphi^{1/r} y> = x^* \varphi^{1/r} <\xi,
\eta>_\mathfrak{H} \varphi^{1/r} y =x^*<\xi, \eta>_\mathfrak{X} y = <\xi x, \eta
y>_\mathfrak{X}.$$
So the correspondence
$$\xi \otimes \varphi^{1/r} x \leftrightarrow \xi x$$
densely defines an $L^p$ module isomorphism
$$\mathfrak{H}_{\mathfrak{X},\varphi} \otimes L^r \simeq \mathfrak{X}.$$ ÊAs
before, the $\mathcal{M}$-action and inner product must agree on the closure,
and the proof is complete.

\bigskip

Since any $L^p$ module is principal, we see that 1) $\| \cdot \|$ is a norm for
$p \ge 1$ and a $p$- (not just $p/2$-) norm for $p <1$; and 2)
$$\bigoplus_c \mathfrak{X}_\alpha = \{ (\xi_\alpha ) | \sum_\alpha <\xi_\alpha,
\xi_\alpha> \in L^{p/2} \};$$
as were mentioned in Section 2.

\bigskip

It also follows from the proof that for any set $S \subset \mathfrak{X}$ an
$L^p$ module,
\begin{equation} \label{E:perp}
S^{\perp \perp} = \overline{S\mathcal{M}}.
\end{equation}
So if $S$ is already an $L^p$ module,
$$\mathfrak{X} = S \oplus_c S^\perp;$$
that is, right $L^p$ submodules are necessarily column summands.

\section{An application to ultraproducts}

Here we give a nontrivial application of Theorem $\ref{T:main}$.

Fix a free ultrafilter $\mathcal{U}$ on $\mathbb{N}$. ÊFor a Banach space $X$,
we define the ultrapower $X_\mathcal{U}$ by
$$\ell^\infty(X)/{\mathcal{N}_\mathcal{U}}, \text{ where }
\mathcal{N}_\mathcal{U} = \{(x_n)| \lim_{n \to \mathcal{U}} \|x_n\| = 0 \}.$$

We will need the following result of Raynaud:

\begin{theorem}
[R] Let $\mathcal{M}$ be a von Neumann algebra. ÊSet
$$\mathcal{N} = Ê((\mathcal{M}_*)_\mathcal{U})^*.$$
Then $\mathcal{N}$ is a von Neumann algebra and
$$L^p(\mathcal{M})_\mathcal{U} \simeq L^p(\mathcal{N}).$$
\end{theorem}
In fact $\mathcal{M}_\mathcal{U}$ is strongly dense in $\mathcal{N}$.

Now take $p \ge 2$ for simplicity, and consider the Banach space
$C^p(\mathcal{M})_\mathcal{U}$. ÊWhile $C^p(\mathcal{M})_\mathcal{U}$ naturally
contains $C^p(\mathcal{N})$, they are not equal; the reader should think of a
sequence of unit vectors in $C^p(\mathcal{M})$ where the support wanders off to
infinity. Ê(Evaluating the limit componentwise is a projection from
$C^p(\mathcal{M})_\mathcal{U}$ onto
$C^p(\mathcal{N}).$) ÊAs huge as $C^p(\mathcal{M})_\mathcal{U}$ is, we can still
gain some control over it via

\begin{proposition}
$C^p(\mathcal{M})_\mathcal{U}$ is a right $L^p$ $\mathcal{N}$-module.
\end{proposition}

\begin{proof}
We explain the $L^p$ module structure. ÊLet $x \in \mathcal{M}_\mathcal{U}
\subset \mathcal{N}$ have representing sequence $(x_n)$, and let $\xi, \eta \in
C^p(\mathcal{M})_\mathcal{U}$ have representing sequences
$(\oplus_{c,k=1}^\infty \xi_k^n),$ $(\oplus_{c,k=1}^\infty \eta_k^n).$ ÊWe
naturally define $\xi x$ by the representing sequence $(\oplus_c \xi_k^n x_n)$;
it is easy to see that this does not depend on the initial choices. ÊSimilarly,
we set
$$<\xi, \eta> = \lim_{n \to \mathcal{U}} < \oplus_{c,k=1}^\infty \xi_k^n,
\oplus_{c,k=1}^\infty \eta_k^n> \in L^{p/2}(\mathcal{M})_\mathcal{U} =
L^{p/2}(\mathcal{N}).$$
It is clear that this inner product generates the norm and is compatible with
the module action.

Finally, we show that the module action extends naturally from
$\mathcal{M}_\mathcal{U}$ to $\mathcal{N}$. ÊThe strong topology that
$\mathcal{M}_\mathcal{U}$ inherits from its action on
$C^p(\mathcal{M})_\mathcal{U}$ is generated by seminorms of the form
\begin{align*}
(x_n) = x &\mapsto \| \xi x \| = \lim_{n \to \mathcal{U}} \| x_n^* <\oplus_c
\xi_k^n, \oplus_c \xi_k^n> x_n \|^{1/2} \\ &= \lim_{n \to \mathcal{U}} \| x_n^*
\sum_k < \xi_k^n, \xi_k^n> x_n \|^{1/2} = \lim_{n \to \mathcal{U}} \|(\sum_k <
\xi_k^n, \xi_k^n> )^{1/2} x_n \| \\ &= \lim_{n \to \mathcal{U}} \| \varphi_n^{1/p} x_n \| = \|\varphi^{1/p} x\|,
\end{align*}
where the $\varphi_n^{1/p} \in L^p(\mathcal{M})$ form a representing sequence for $\varphi^{1/p} \in L^p(\mathcal{N}).$ ÊBy Lemma \ref{T:strong}, these are exactly the seminorms which generate the strong topology on
$\mathcal{M}_\mathcal{U}$ inside $\mathcal{N}$. ÊThis completes the proof.
\end{proof}

By Theorem $\ref{T:main}$, we know that any $L^p$ module can be written as a
column sum. ÊOne can think of $C^p(\mathcal{M})_\mathcal{U}$ as containing
countably many copies of $L^p(\mathcal{N})$ from componentwise limits, plus
uncountably many more from all the directions in which support might wander.
Perhaps it is more natural to think of $C^p(\mathcal{M})_\mathcal{U}$ as a
continuous column integral of $L^p(\mathcal{N})$ over a \textit{very} large
space; the adventurous reader may want to consider how to make this statement
more precise.

\section{Commutants and categorical properties}

Consider a countably generated right $L^p$ module $\mathfrak{X}$ for a $\sigma$-finite von Neumann algebra $\mathcal{M}$. ÊBy Theorem
$\ref{T:main}$, there are projections $\{q_n \}$ with
$$\mathfrak{X} \simeq \bigoplus_c q_n L^p \simeq (\sum q_n \otimes e_{nn}) C^p = qC^p.$$ ÊIf two such projections $q_1, q_2 \in M_\infty(\mathcal{M})$ are Murray-von Neumann equivalent via a partial isometry $v$, we have a module isomorphism (even isometric) $q_1 C^p \simeq q_2 C^p$ via left multiplication by $v$. ÊWe will obtain the converse after proving

\begin{proposition} \label{T:pcomm}
On $q C^p,$ the right action of $\mathcal{M}$ ($= R(\mathcal{M})$) and
the left action of $qM_\infty(\mathcal{M})q$ ($= L(qM_\infty(\mathcal{M})q)$) are commutants of each other.
\end{proposition}

\begin{proof}

Once again we may assume that the module is faithful for the right action of $\mathcal{M}$. ÊNow note that the two actions mentioned are commuting and bounded: boundedness of $L(qM_\infty(\mathcal{M})q)$ follows by viewing it as a subalgebra of $M_\infty(\mathcal{M})$ and using the H\"{o}lder inequality. ÊFinally, by the remark before the proposition, we may assume that $q^\perp \sim 1$.

So let $T$ be a bounded operator on $q C^p$ commuting with the right action of $\mathcal{M}$, and then set $T' = T \circ L(q).$ Ê$T'$ is bounded and commutes with $R(\mathcal{M})$ on all of $C^p$.

Since $T'$ acts on column vectors, it has a matrix representation as $(T'_{ij})$, where each $T'_{ij}$ operates on $L^p(\mathcal{M}).$ ÊFix a single $T'_{ij}$. ÊFor any $\xi \in L^p(\mathcal{M})$ and $x \in \mathcal{M}$ we may consider the vector $\bar{\xi}$ in $C^p$ with $\xi$ in the $j$th position and 0 elsewhere. ÊSince $T'$ commutes with $R(\mathcal{M})$,
$$(T'_{ij} \xi) x = ((T' \bar{\xi}) x)_i = (T'(\bar{\xi} x))_i = T'_{ij}(\xi x).$$
By Corollary \ref{T:comm} we know that $T'_{ij} = L(y_{ij})$ for some $y_{ij} \in \mathcal{M}$. ÊConsidering the kernel and range, we deduce that $y_{ij} \in q_i \mathcal{M} q_j.$ ÊThen $T' = L((y_{ij}))$ for some bounded operator $(y_{ij}) \in q M_\infty(\mathcal{M})q$, and this representation is the restriction $T$ as well.

Instead of trying to check that any operator commuting with $L(qM_\infty(\mathcal{M})q)$ must lie inside $R(\mathcal{M})$, we give a small argument involving projections in order to invoke symmetry. ÊIn the $\sigma$-finite algebra $M_\infty \otimes M_\infty \otimes \mathcal{M}$, the projections
$$I \otimes q \text{ and } e_{11} \otimes I \otimes I_\mathcal{M}$$
are both properly infinite and therefore equivalent ([KR], Corollary 6.3.5). ÊThey remain so after subtracting their common subprojection $e_{11} \otimes q$ (we assumed $q^\perp \sim I \otimes I_\mathcal{M}$), allowing us to find a partial isometry $v$ between them which fixes $e_{11} \otimes q$. ÊConjugation by $v$ gives an isomorphism
\begin{align} \label{E:corners}
M_\infty(q M_\infty(\mathcal{M}) q) &= (I \otimes q)( M_\infty \otimes M_\infty \otimes \mathcal{M}) (I \otimes q ) \\ \notag &\simeq(e_{11} \otimes I \otimes I_\mathcal{M}) (M_\infty \otimes M_\infty \otimes \mathcal{M} ) ( e_{11} \otimes I \otimes I_\mathcal{M}) \\ \notag &= M_\infty(\mathcal{M}).
\end{align}
Now let $r= v^*(e_{11} \otimes e_{11} \otimes 1_\mathcal{M}) v$ be the projection in the first algebra which corresponds to $e_{11} \otimes 1_\mathcal{M}$ in the last, and notice $e_{11} \otimes q$ is the ``outer" matrix unit $e_{11}'$ for $M_\infty(q M_\infty(\mathcal{M}) q).$
Via the isomorphism above, we have the isomorphic bimodule presentations
\begin{equation} \label{E:basealg}
q(M_\infty(\mathcal{M}))q \:- \: q L^p(M_\infty(\mathcal{M}))e_{11} \: - \:\mathcal{M} \quad
\simeq
\end{equation}
$$q(M_\infty(\mathcal{M}))q \: - \: e_{11}' L^p(M_\infty(q(M_\infty(\mathcal{M}))q)) r\: -
\:r(M_\infty(q(M_\infty(\mathcal{M}))q))r.$$
(The point is to observe that module and commutant are written as reduced
amplifications of the \textit{left} algebra.) ÊNow applying the first argument
finishes the proof.
\end{proof}

That $q$ be diagonal, i.e. of the form $\sum q_n \otimes e_{nn}$, is actually unnecessary. ÊFor any projection $q$ in $M_\infty (\mathcal{M})$, the $L^p$ module $q C^p$ inherits its structure from $C^p.$

\begin{corollary}
If the $L^p$ modules $q_1 C^p$, $q_2 C^p$ are isomorphic, then the projections $q_1$, $q_2$ are Murray von-Neumann equivalent.
\end{corollary}

\begin{proof}
The proof is no different than the $L^2$ case. ÊIf $S$ is the isomorphism, extend it to
$$\widetilde{S}: q_1 C^p \oplus_c q_1^\perp C^p = C^p \to C^p;$$ $$\widetilde{S}(\xi \oplus_c \eta) = S(\xi).$$
$\widetilde{S}$ is clearly bounded, so by Proposition \ref{T:pcomm}, it is given by left composition with some $y \in M_\infty(\mathcal{M})$. ÊBy considering the kernel and range of $\widetilde{S}$, we see that $s_r(y) = q_1$ and $s_\ell(y) = q_2.$
\end{proof}

By virtually the same argument we obtain

\begin{corollary} \label{T:geninter}
If $\frac{1}{p} + \frac{1}{r} = \frac{1}{q}$,
$$\text{Hom}(q_1 C^p_\mathcal{M}, q_2 C^q_\mathcal{M}) = L(q_2 L^r(M_\infty(\mathcal{M})) q_1).$$
\end{corollary}

\bigskip

\textit{Remark:} Proposition \ref{T:pcomm} and its corollaries still hold without the assumptions that the algebra is $\sigma$-finite and the module is countably generated. Ê(This requires either a direct limit argument or a more subtle calculation with projections.) ÊIn the general case, the typical module is $q L^p(M_J(\mathcal{M})) e_{11}$ for some cardinal $J$ and projection $q \in M_J(\mathcal{M});$ we will not need the full result in the sequel and so opted for clarity.

\bigskip

Now we continue to investigate the category of isomorphism classes of countably generated right $L^p$ $\mathcal{M}$-modules, with intertwiners as morphisms, which we call $Right\: L^p Mod(\mathcal{M}).$ ÊThese are submodules of $C^p$; from the foregoing discussion we may conclude that they are parameterized by Murray von-Neumann equivalence classes of projections in $M_\infty(\mathcal{M})$, which is $V(M_\infty(\mathcal{M}))$ in the language of $K$-theory [W-O]. ÊIt should be clear
that this is an additive category, with addition being the column sum of orthogonal representatives. ÊThis actually gives us monoidal equivalence with $V(M_\infty(\mathcal{M}))$:
$$q_1 C^p \oplus_c q_2 C^p = (q_1 + q_2) C^p, \quad q_1 \perp q_2$$
corresponds exactly to
$$[q_1] + [q_2] = [q_1 + q_2], \quad q_1 \perp q_2.$$
It follows that $C^p \oplus_c qC^p \simeq C^p$, which is the $L^p$ version of Kasparov's stabilization theorem for Hilbert C*-modules ([L], Theorem 6.2). ÊIn case $\mathcal{M}$ is a $\text{II}_1$
factor, we can make the correspondence with $V(M_\infty(\mathcal{M})) \simeq [0,\infty]$ explicit with the natural definition
\begin{equation}
\text{dim}_\mathcal{M} \: (q C^p) = \tau_{M_\infty(\mathcal{M})} (q).
\end{equation}
Clearly $\text{dim}_{\mathcal{M}}(\oplus_c \mathfrak{X}_i) = \sum
\text{dim}_{\mathcal{M}} \mathfrak{X}_i$.

All of this is identical to the $L^2$ case, but we recall the difference at the
vector level: the norm in a column sum ($p \ne 2$)
$$\|\xi \oplus_c \eta Ê\| = \|\xi^* \xi + \eta^* \eta \|^{1/2}$$
is not, in general, a function of the norms in each component. ÊNow it may occur
to the reader to try a ``diagonal" sum
$\left(
\begin{smallmatrix}
\mathfrak{X} & 0 \\
0 & \mathfrak{Y}
\end{smallmatrix}
\right)$, as is done for operator spaces. ÊThis \textit{is} an $\ell^p$ direct
sum, but no compatibility is required or retained: the diagonal sum of a right
$L^p$ $\mathcal{M}_1$-module and a right $L^p$ $\mathcal{M}_2$-module is a right
$L^p$ $(\mathcal{M}_1 \oplus \mathcal{M}_2)$-module. ÊIf $\mathcal{M}_1 =
\mathcal{M}_2 = \mathcal{M}$, the diagonal sum is algebraically an
$\mathcal{M}$-module, but not necessarily an $L^p$ module in our sense: the
inner product would naturally be $L^{p/2}(\mathcal{M} \oplus
\mathcal{M})$-valued. ÊThe difference is already apparent in the simplest
possible case: $$\mathcal{M} = \mathfrak{X} = \mathfrak{Y} = \mathbb{C}.$$ ÊAs
modules,
\begin{equation}
\mathbb{C} \oplus_c \mathbb{C} =
\left(
\begin{smallmatrix}
\mathbb{C} \\
\mathbb{C}
\end{smallmatrix} \right) {}_\mathbb{C}, \quad
\left\| \left(
\begin{smallmatrix}
a \\
b
\end{smallmatrix}
\right) \right\| = (|a|^2 + |b|^2)^{1/2},
\end{equation}
$$\mathbb{C} \oplus_d \mathbb{C} = \left(
\begin{smallmatrix}
\mathbb{C} & 0 \\
0 & \mathbb{C}
\end{smallmatrix} \right) {}_{\mathbb{C}I_2 \subset M_2}, \quad
\left\| \left(
\begin{smallmatrix}
a & 0 \\
0 & b
\end{smallmatrix}
\right) \right\| = (|a|^p + |b|^p)^{1/p}.$$
$\mathbb{C} \oplus_d \mathbb{C}$ cannot be a right $L^p$ $\mathbb{C}$-module,
since it is apparently not isometrically isomorphic to the only two-dimensional
right $L^p$ $\mathbb{C}$-module, $\mathbb{C} \oplus_c \mathbb{C}$. Ê(Instead it
is a right $L^p$ $(\mathbb{C} \oplus \mathbb{C})$-module.)

\bigskip

Since there are many equivalent constructions of $L^p(\mathcal{M})$, it should
not be surprising that there are other ways to build the class of $L^p$ modules. ÊWe do not reproduce the details from [S1] but simply note that the class of
countably generated right $L^p$ modules, modulo spatial isomorphism, can also be
described as
\begin{itemize}
\item a minimal class of complete right $\mathcal{M}$-modules which contains
$L^p(\mathcal{M})$ and is closed under taking submodules and forming countable
column sums (recall equation $\eqref{E:perp}$);

\item a class of spaces of ``column" operators which satisfy a $-1/p$-homogeneity condition in the sense of Connes-Hilsum [Hi];

\item a class of interpolation spaces, following Kosaki [K1].
\end{itemize}

\bigskip

In the sequel we will frequently be concerned with left actions. ÊOf course, the theory of left $L^p$ modules is entirely analogous. ÊThe counterparts to column sums, $C^p$, and $Right \: L^p Mod$ we call \textit{row sums}, $R^p$, and $Left \: L^p Mod$. ÊThere is a 1-1 correspondence between left and right $L^p$ $\mathcal{M}$-modules given by the \textit{contragredient} ${}_\mathcal{M} \bar{\mathfrak{X}}$ of $\mathfrak{X}_\mathcal{M}$: $\bar{\mathfrak{X}}$ is conjugate linearly isomorphic to $\mathfrak{X}$, with left action $x \cdot \bar{\xi} = \overline{\xi x^*}$ and inner product $<\bar{\xi}, \bar{\eta}> = <\xi, \eta>$. ÊOf course, one may similarly take the contragredient of a left $L^p$ module; $\bar{\bar{\mathfrak{X}}}$ is canonically isomorphic to $\mathfrak{X}$. ÊIt is easy to see that when $\mathfrak{X}$ is represented as a principal $L^p$ module, the contragredient corresponds to the operator adjoint.

\begin{definition}
An \textbf{$L^p \: \mathcal{M}-\mathcal{N}$ bimodule} is an
$\mathcal{M}-\mathcal{N}$ bimodule (meaning that the actions commute) which is
simultaneously a left $L^p$ $\mathcal{M}$-module and a right $L^p$
$\mathcal{N}$-module. ÊWe denote the category of isomorphism classes, with intertwiners as morphisms, by $L^p Bimod(\mathcal{M}, \mathcal{N}).$
\end{definition}

Notice, by Proposition $\ref{T:pcomm}$ and $\eqref{E:basealg}$, that every left
or right $L^p$ $\mathcal{M}$-module is an $L^p$ bimodule, with opposite action
coming from the commutant. ÊWe will explore this more fully in the next section.

Our final observation of this section concerns the relative tensor product, a
sort of ``multiplication" for Hilbert modules. ÊThe original
arguments are due to Connes and Sauvageot (and found in [P] and [Sa]); the informed reader will recognize our $L^p$ formulation as a minor modification. ÊAs explained in [S3], on the module level the relative tensor product only ``sees" the projections (more precisely, the elements of $V(M_\infty(\mathcal{M})))$ which determine the modules. ÊThe densities of the modules - all $1/2$ in the usual case - are irrelevant, and so we may choose any $p,q,r$ we please. ÊIn the following definition the notations $\mathcal{L}(\mathfrak{X}_\mathcal{M})$ and $\mathcal{L}({}_\mathcal{M} \mathfrak{Y})$ stand for commutants.

\begin{definition}
By an \textbf{$(\mathcal{M},p,q,r)$-relative tensor product} we mean a functor,
covariant in both variables,
$$Right\:L^p Mod(\mathcal{M}) \times Left\:L^q Mod(\mathcal{M}) \to L^r Bimod(\mathcal{L}(\mathfrak{X}_\mathcal{M}),\mathcal{L}({}_\mathcal{M}
\mathfrak{Y}));$$
$$(\mathfrak{X},\mathfrak{Y}) \mapsto \mathfrak{X} \otimes_{\mathcal{M},p,q,r}
\mathfrak{Y}$$
which satisfies
$$L^p(\mathcal{M}) \otimes_{\mathcal{M},p,q,r} L^q(\mathcal{M}) \simeq
L^r(\mathcal{M})$$
as bimodules.
\end{definition}

We remind the reader that the sums in these categories (which the relative
tensor product must distribute, by functoriality) are not direct. ÊSo, for
example,
$$(\oplus_c L^p) \otimes_{\mathcal{M},p,q,r} L^q \simeq \oplus_c(L^r).$$
By decomposition and functoriality, it is simple to see that such functors exist
and are unique up to unitary equivalence. ÊOne has the following representation
result:
\begin{proposition} \label{T:colxrow}
Let $\mathfrak{X} \simeq q_1 C^p \in Right \:L^p Mod(\mathcal{M}),\:
\mathfrak{Y} \simeq ÊR^q q_2 \in Left\:L^q Mod(\mathcal{M})$
for some $q_1,q_2 \in \mathcal{P}(M_\infty(\mathcal{M})).$ ÊThen
$$\mathfrak{X} \otimes_{\mathcal{M},p,q,r} \mathfrak{Y} \simeq
q_1L^r(M_\infty(\mathcal{M}))q_2$$
with natural action of the commutants.
\end{proposition}

It is also possible to give an element-wise construction of the relative tensor
product based on a fixed faithful state (or weight) $\varphi$. ÊThe usual
construction is
$$\xi \otimes_\varphi \eta = \xi \varphi^{-1/2} \eta$$
for a suitable dense set of $\xi, \eta$, and the $(p,q,r)$-relative tensor
product requires
$$\xi \otimes_\varphi \eta = \xi \varphi^{\frac{1}{r} - \frac{1}{p} -
\frac{1}{q}} \eta.$$
Both of these identities are discussed in [S2], and in [S3] the preclosedness of this relative tensor map is investigated in full. ÊThe reader will notice that the auxiliary constructions introduced to prove the $p<2$ case of Theorem $\ref{T:main}$ are nothing but relative tensor products.

\section{$L^p$ bimodules}

The theory of $L^2$ bimodules, which contains that of subfactors, is one of the
most fruitful fields in the study of von Neumann algebras. ÊBut for the $L^p$ analogues with $p \ne 2$, the lack of Hilbert space symmetry makes for a much more restrictive theory. ÊOne deficit which is apparent from the
outset is that $L^p$ bimodules do not add: row and column sums preserve one
algebra only. ÊWe will see that there are other significant limitations.
\textit{In this section we simplify the discussion by assuming $1 < p<\infty$, $p\ne 2$, all algebras to be $\sigma$-finite, and all $L^p$ modules to be countably generated and faithful.}

\bigskip

The structure theorems proven so far show that every left or right $L^p$
$\mathcal{M}$-module is an $L^p$ bimodule, with opposite action coming from the
commutant. ÊBut of course the commutant is the ``largest" choice, so an
$\mathcal{M}-\mathcal{N}$ $L^p$ bimodule gives injective homomorphisms of each
algebra into an amplification of the other.

\begin{lemma} \label{T:incl}
Let $\mathfrak{X}$ be an $\mathcal{M}-\mathcal{N}$ $L^p$ bimodule, and suppose
that $\mathfrak{X} \simeq R^p(\mathcal{M}) q$ via an isomorphism $T$. ÊThen
there is an injective normal *-homomorphism
$$\pi: \mathcal{N} \hookrightarrow qM_\infty(\mathcal{M})q$$
such that
$$T(\xi n) = T(\xi) \pi(n).$$
\end{lemma}

\begin{proof}
The only things to show are that $\pi$ is normal and *-preserving. ÊFor
normality, note that by Lemma $\ref{T:strong}$ the strong topologies on the unit
balls of $\mathcal{N}$ and $\pi(\mathcal{N})$ are both given by the module
actions, which are identical.

A variant of Lemma \ref{T:spec} says that the norms are also generated by the
module actions. ÊSince (orthogonal) projections are exactly idempotents of norm one, $\pi$ takes projections to projections. ÊBy approximating with projections,
we see $\pi$ takes self-adjoint elements to self-adjoint elements. ÊFinally,
$$\pi(x^*) = \pi(\text{Re }x) - i\pi(\text{Im }x) = \pi(x)^*.$$
\end{proof}

One is tempted to follow the theory of correspondences [P] and guess that $L^p$
bimodules are equivalent to normal unital *-homomorphisms, but this is asking
too much. ÊFor example, $L^p(\mathcal{M})$ is an $\mathcal{M}-\mathcal{M}$ $L^p$
bimodule, and naturally $\mathbb{C} \subset \mathcal{M}$. ÊBut if
$L^p(\mathcal{M})$ is a right $L^p$ $\mathbb{C}$-module, then $L^p(\mathcal{M})
\simeq q C^p(\mathbb{C})$ for some $q$, and this last object is actually a
Hilbert space. ÊThen the norm in $L^p(\mathcal{M})$ follows the parallelogram
law; if $\varphi$, $\psi$ are states with orthogonal supports,
\begin{equation} \label{E:para}
4 = 2(\|\varphi^{1/p}\|^2 +\|\psi^{1/p} \|^2) = \| \varphi^{1/p} + \psi^{1/p}
\|^2 + \|\varphi^{1/p} - \psi^{1/p} \|^2 = 2 \cdot 2^{2/p}.
\end{equation}
Since this is false when $p \ne 2$, such states cannot exist. ÊThus
$L^p(\mathcal{M})$ is a right $L^p$ $\mathbb{C}$-module iff $\mathcal{M} =
\mathbb{C}$.

The cause of such a phenomenon is clear: noncommutative $L^p$ spaces remember
their generating algebras (except for amnesiac $p=2$). ÊThe existence of an
$L^p$ $\mathcal{M}-\mathcal{N}$ bimodule therefore implies a relationship
between $\mathcal{M}$ and $\mathcal{N}$, and the remainder of this section is
devoted to a precise description of this relationship.

A major tool will be the following result of Raynaud and Xu (relying heavily on Kosaki's papers [K2], [K3], where a subcase is proved).

\begin{theorem} \label{T:clarkson} [RX]
When $p \ne 2$, two elements $\xi, \eta$ of a noncommutative $L^p$ space satisfy $s_\ell(\xi) \perp
s_\ell(\eta), s_r(\xi) \perp s_r(\eta)$ iff they satisfy
\begin{equation} \label{E:clarkson}
\|\xi + \eta\|^p + \|\xi - \eta\|^p = 2\left( \|\xi\|^p + \|\eta\|^p \right).
\end{equation}
\end{theorem}

\begin{proposition} \label{T:center}
Let $\mathfrak{X} \simeq R^p(\mathcal{M}) q$ be an $L^p$
$\mathcal{M}-\mathcal{N}$ bimodule. ÊThen the centers of $\mathcal{M}$ and
$\mathcal{N}$ are isomorphic and act identically on $\mathfrak{X}$.
\end{proposition}

\begin{proof}
We will use the classical $L^p$ notion of \textit{$L^p$-projection} [B]: an
idempotent $E$ on a Banach space satisfying
\begin{equation} \label{E:lpproj}
\|\xi\|^p = \|E\xi\|^p + \|\xi - E\xi\|^p
\end{equation}
for all elements $\xi.$ ÊCommutative $L^p$ spaces are characterized by having
sufficiently many $L^p$-projections, but the same is not true in the
noncommutative setting. ÊWe will show that the $L^p$-projections on
$\mathfrak{X}$ can be identified spatially with the central projections of
either $\mathcal{M}$ or $\mathcal{N}$.

That a central projection is an $L^p$-projection is clear. ÊSo let $E$ be an $L^p$-projection, and make the identification $\mathfrak{X} \simeq
R^p(\mathcal{M}) q$. ÊFor any two vectors $\xi \in \text{Ran}E, \, \eta \in \ker E$, \eqref{E:lpproj} implies \eqref{E:clarkson} (since $-\eta \in \ker E$ also).
ÊThen $\xi$ and $\eta$ have orthogonal left and right supports, and we have the
decomposition in $\mathcal{M}$
$$\bigvee_{\xi \in \text{Ran}E} s_\ell(\xi) + \bigvee_{\eta \in \ker E}
s_\ell(\eta) = 1.$$
Now the left action of the first projection above is apparently $E$, which must
also equal the right action of $\bigvee_{\xi \in \text{Ran}E} s_r(\xi) \in
(\mathcal{M}')^{\text{op}}= qM_\infty(\mathcal{M})q.$ ÊThis is only possible if
the projection in $\mathcal{M}$ is central.

Thus each central projection in $\mathcal{N}$, being an $L^p$-projection, is
identified spatially with a central projection in $\mathcal{M}$. ÊIt follows that the centers of $\mathcal{M}$ and $\mathcal{N}$ are isomorphic.

\end{proof}

At this point, our original approach was to decompose $\mathfrak{X}$ into a
direct integral of $L^p$ bimodules between factors. ÊThis requires a significant
detour into measure theory, and we have opted to omit these arguments (which may
appear elsewhere) and deal with central projections.

We will use the following two results, which characterize isometries between noncommutative $L^p$ spaces under certain circumstances.

\begin{theorem}
[Ye] \label{T:yeadon} Let $\{\mathcal{M}, \tau_{\mathcal{M}} \}, \{\mathcal{N},
\tau_{\mathcal{N}} \}$ be semifinite von Neumann algebras with given traces.
Let $T$ be a linear isometry from $L^p(\mathcal{M}, \tau_{\mathcal{M}})$ to
$L^p(\mathcal{N}, \tau_{\mathcal{N}})$ $(p \ne 2)$, where we view these as spaces of
$\tau_{\mathcal{M}}$ or $\tau_{\mathcal{N}}$-measurable operators. ÊThen there
exist, uniquely, a partial isometry $w \in \mathcal{N}$, an injective normal Jordan
*-homomorphism $J$ of $\mathcal{M}$ into $\mathcal{N}$, and a positive unbounded
operator $B$ affiliated with $J(\mathcal{M})' \cap \mathcal{N} $, all satisfying
$$w^*w = J(1) = s(B);$$
$$\tau_{\mathcal{M}}(x) = \tau_{\mathcal{N}}(B^p J(x)), \qquad \forall x \in
\mathcal{M}_+;$$
$$T(x) = wB J(x), \qquad \forall x \in \mathcal{M} \cap L^p(\mathcal{M},
\tau_{\mathcal{M}}).$$
\end{theorem}

It is worth explaining here that a map between von Neumann algebras is \textit{Jordan} if it preserves the Jordan product $x \circ y = (1/2)(xy + yx)$. ÊAn injective normal Jordan *-homomorphism is the sum of a *-homomorphism and a *-antihomomorphism; the two supports, which are central projections, have sum $\geq 1$, and the two ranges are orthogonal. ÊA surjective Jordan *-isomorphism is necessarily normal and so can be centrally decomposed, in both domain and range, into a *-isomorphism and a *-antiisomorphism. Ê(See [HaS] for details. ÊThe only ambiguity in these decompositions arises from abelian summands.)

\begin{theorem} [S4] \label{T:ncbst}
Suppose that $T: L^p(\mathcal{M}) \to L^p(\mathcal{N})$ is a surjective isometry, $1 < p < \infty$, $p \ne 2$. ÊThen there are a surjective Jordan *-isomorphism $J: \mathcal{M} \to \mathcal{N}$ and a unitary $u \in \mathcal{N}$ such that
\begin{equation} \label{E:ncbst}
T(\varphi^{1/p}) = u(\varphi \circ J^{-1})^{1/p}, \qquad \forall \varphi \in \mathcal{M}_*^+.
\end{equation}
We have $T(\xi x) = T(\xi) J(x)$ (resp. $T(\xi x) = u J(x) u^* T(\xi)$) when $\xi x$ is supported on a central summand for which $J$ is multiplicative (resp. antimultiplicative).
\end{theorem}

\begin{proposition} \label{T:jordan}
Let $\mathfrak{X} \simeq R^p(\mathcal{M}) q$ be an $L^p$
$\mathcal{M}-\mathcal{N}$ bimodule with $q \preccurlyeq e_{11}$. ÊThen the inclusion
$$\pi:\mathcal{N} \hookrightarrow q\mathcal{M}q$$
is surjective.
\end{proposition}

\begin{proof}
If necessary, implement an isomorphism so that $q \le e_{11}.$ ÊThe hypotheses
mean that we have bimodules
\begin{equation} \label{E:tab1}
\mathcal{M} - L^p(\mathcal{M})q - q\mathcal{M}q,
\end{equation}
$$q'M_\infty(\mathcal{N})q' - q' C^p(\mathcal{N}) - \mathcal{N},$$
where $q' \in \mathcal{P}(M_\infty({N}))$ and the bimodules are isometrically
isomorphic as Banach spaces. ÊBy Lemma $\ref{T:incl}$ we have inclusions
consistent with the module actions:
\begin{equation} \label{E:incl}
\mathcal{M} \subset q'M_\infty(\mathcal{N})q', \quad \mathcal{N}
\overset{\pi}{\subset} q\mathcal{M}q.
\end{equation}
So via its left action, the projection $q \in \mathcal{M}$ is identified with a projection $q'' \in
q'M_\infty(\mathcal{N})q'$. ÊIf we implement $L(q) \leftrightarrow L(q'')$, we get subbimodules
\begin{equation} \label{E:tab2}
q\mathcal{M}q\: - \: L^p(q\mathcal{M}q) \: - \: q\mathcal{M}q,
\end{equation}
$$q''M_\infty(\mathcal{N})q'' \: - \: q'' C^p(\mathcal{N}) \: - \: \mathcal{N}.$$
Again the modules themselves are isometrically isomorphic, and we still have the same inclusion $\mathcal{N} \overset{\pi}{\subset} q\mathcal{M}q$ from the right actions.

Our next step will be to show that $q'' \sim e_{11}^\mathcal{N}$, so that $q'' C^p(\mathcal{N})$ may be replaced with $L^p(\mathcal{N})$ in \eqref{E:tab2}. ÊBy implementing projections from the common center of $\mathcal{M}$ and $\mathcal{N}$, we may consider separately the cases where $\mathcal{N}$ is finite or properly infinite. ÊIf $\mathcal{N}$ is properly infinite, then the given center-preserving inclusions $\mathcal{N} \hookrightarrow q\mathcal{M}q \hookrightarrow q''M_\infty(\mathcal{N})q''$ imply that $q''$ is properly infinite, and so $q'' \sim e_{11}^\mathcal{N}$.

If $\mathcal{N}$ is finite, first we argue that $q\mathcal{M}q$
must be finite. ÊFor otherwise $L^p(q\mathcal{M}q)$ contains an
isometric copy of the Schatten class $S^p (=
L^p(\mathcal{B}(\mathfrak{H}), \text{tr}))$, and so we have that $S^p$
embeds isometrically in $q''C^p(\mathcal{N})$, where
$\mathcal{N}$ is finite. ÊLetting $\tau$ be a faithful normal trace on $\mathcal{N}$, we show that this is impossible.  Indeed, if
$p>2$, it is easily seen that $C^p(\mathcal{N})$ is a subspace of
the intersection Ê$X=L^p(M_\infty(\mathcal{N}), \tau \otimes \text{tr}) \cap
L^2(M_\infty(\mathcal{N}), \tau \otimes \text{tr})$.  According to \cite{J2}, $X$ embeds into $L^p(\widetilde{\mathcal{M}})$ for some finite $\widetilde{\mathcal M}$.  This yields an embedding of $S^p$ into $L^p(\widetilde{\mathcal M})$, which is absurd in view of the result of
Sukochev \cite[Theorem 3.1]{Su}.  For $1 < p<2$, we replace $X$ by the sum $L^p(M_\infty(\mathcal{N}), \tau \otimes \text{tr}) + L^2(M_\infty(\mathcal{N}), \tau \otimes \text{tr})$.  Again by \cite{J2}
$X$ embeds into $L^p(\widetilde{\mathcal M})$ for some finite
$\widetilde{\mathcal M }$ and thus an embedding of $S^p$ in
$C^p({\mathcal N})$ provides an embedding of $S^p$ into
$L^p(\widetilde{\mathcal M})$.  In this case, we may refer to the main result of \cite{HRS} for the fact that this is impossible.  So $q\mathcal{M}q$ is finite with faithful normal trace $\tau'$, and we may apply Theorem \ref{T:yeadon} to the $L^p$ isometry
$$T: L^p(q\mathcal{M}q, \tau') \simeq q'' C^p(\mathcal{N}) \hookrightarrow
L^p(M_\infty(\mathcal{N}), \tau \otimes \text{tr}).$$
With $T = w B J(\cdot)$, the conditions of the theorem imply that $q'' = s_\ell(w) \sim s_r(w) = e_{11}^\mathcal{N}$, as desired.

This means that we may replace the bottom line of \eqref{E:tab2} by $\mathcal{N} - L^p(\mathcal{N}) - \mathcal{N}$. ÊWe still have $\pi: \mathcal{N} \hookrightarrow q\mathcal{M}q$, and we set $S: L^p(\mathcal{N}) \to L^p(q\mathcal{M}q)$ to be the isometric isomorphism of Banach spaces. ÊApplying Theorem \ref{T:ncbst} to $S$ we find the underlying pair $u \in q\mathcal{M}q, J: \mathcal{N} \to q\mathcal{M}q$; let $z, z^\perp$ be central projections of $\mathcal{N}$ which divide $J$ into multiplicative and antimultiplicative parts. ÊThe intertwining relation between $S$ and $\pi$ gives
\begin{align} \label{E:twining}
u (\varphi &\circ J^{-1})^{1/p} \pi(x) = S(\varphi^{1/p}) \pi(x) = S(\varphi^{1/p} x) \\
\notag &= u (\varphi \circ J^{-1})^{1/p} J(xz) + u J(xz^\perp)(\varphi \circ J^{-1})^{1/p}, \qquad \varphi \in \mathcal{N}_*^+, \: x \in \mathcal{N}.
\end{align}
We observe that $J$ and $\pi$ both identify the centers, so we may multiply on the left by $\pi(z^\perp) u^*$ to get
$$(\varphi \circ J^{-1})^{1/p} \pi(xz^\perp) = ÊJ(xz^\perp)(\varphi \circ J^{-1})^{1/p}, \qquad \varphi \in \mathcal{N}_*^+, \: x \in \mathcal{N}.$$
Then $R(\pi(xz^\perp)) = L(J(x z^\perp))$ is central for any $x$, so that $\mathcal{N}z^\perp$ is abelian, and in fact $J$ is multiplicative. ÊNow \eqref{E:twining} shows that $\pi = J$. ÊSince we know that $J$ is surjective, this finishes the proof.

\end{proof}

Return now to the general situation of an $L^p$ $\mathcal{M}$-$\mathcal{N}$
bimodule $\mathfrak{X} \simeq R^p(\mathcal{M}) q \simeq q'C^p(\mathcal{N}),$ and identify the centers of $\mathcal{M}$ and $\mathcal{N}$. ÊUse the comparability theorem to
find the largest central projections $z,z'$ satisfying $zq \preccurlyeq
ze_{11}^\mathcal{M}, z'q' \preccurlyeq z'e_{11}^\mathcal{N}$. ÊWith $z'' = z \vee
z'$, Proposition $\ref{T:jordan}$ tells us that $z'' \mathcal{M}$ and $(z''
\mathcal{N})^{\text{op}}$ are commutants on $z''\mathfrak{X}.$ ÊOn every central
summand of the complement, both $p$ and $q$ are strictly larger than
$e_{11}^\mathcal{M}$ and $e_{11}^\mathcal{N},$ respectively. ÊIt follows that $z''^\perp \mathcal{M}$ and $z''^\perp \mathcal{N}$ are finite; we will show that in fact they are abelian.

\begin{proposition}
Let $\mathcal{M}$ and $\mathcal{N}$ be finite algebras, and assume that
$\mathcal{M}$ has no abelian central summand. ÊIf $q \in
\mathcal{P}(M_\infty(\mathcal{M}))$ and $q' \in
\mathcal{P}(M_\infty(\mathcal{N}))$ are projections such that $qz \succneqq
e_{11}^\mathcal{M}z$ and $q'z' \succneqq e_{11}^\mathcal{N}z'$ for all central projections $z \in M_\infty(\mathcal{M})$, $z' \in M_\infty(\mathcal{N})$, then there is no $\mathcal{M}-\mathcal{N}$ $L^p$ bimodule $\mathfrak{X} \simeq R^p(\mathcal{M})q \simeq q'C^p(\mathcal{N}).$
\end{proposition}

\begin{proof}
Seeking a contradiction, let $\mathfrak{X}$ be such a bimodule. ÊChoose finite traces $\tau_\mathcal{M}, \tau_\mathcal{N}$ and consider the $L^p$ elements to be measurable operators. ÊAs before, we may assume that $q \ge e_{11}^\mathcal{M}$ and $q' \ge e_{11}^\mathcal{N}.$ ÊLet $T$ be the isometric isomorphism from
$R^p(\mathcal{M})q$ to $q'C^p(\mathcal{N})$; the domain naturally contains $R^p
e_{11}^\mathcal{M} \simeq L^p(\mathcal{M})$ to give the isometric restriction
$$T_1: L^p(\mathcal{M}, \tau_\mathcal{M}) \cap \mathcal{M} \ni x \mapsto T((\begin{smallmatrix} x & 0 & \cdots
\end{smallmatrix})) \in q'C^p(\mathcal{N}) \subset L^p(M_\infty(\mathcal{N}), \tau_\mathcal{N} \otimes \text{tr}).$$

The vector $\xi = (\begin{smallmatrix} 1 & 0 & \cdots \end{smallmatrix}) \in R^p(\mathcal{M})q$ has full left support, so by Theorem \ref{T:clarkson} it satisfies equation \eqref{E:clarkson} for no nonzero $\eta \in \mathfrak{X}$. ÊSince $T$ is an isometric identification, the same is true for $T(\xi) \in q'C^p(\mathcal{N})$. ÊNow our assumption on the size of $q'$ means that $T(\xi)$ can have full left support on no central summand, so by Theorem \ref{T:clarkson} again we must have $s_r(T(\xi)) = e_{11}^\mathcal{N} = 1_\mathcal{N}$.

Theorem \ref{T:yeadon} tells us that $T_1(x) = w B J(x).$ Ê\textit{A priori} these operators are affiliated with $M_\infty(\mathcal{N})$, but the conditions in the theorem imply
$$w B = T_1 (1) \in L^p(M_\infty(\mathcal{N}))e_{11} \Rightarrow J(1) = s_r(w) =
s(B) = e_{11}^\mathcal{N}.$$
Thus we see that $B$ is affiliated with, and $J(\mathcal{M})$ are elements in, $e_{11}M_\infty(\mathcal{N})e_{11}$. ÊWe naturally identify the latter algebra with $\mathcal{N}$, so that $J$ is unital.

Choose $y \in R^p(\mathcal{M})(q-e_{11})$ with $s_\ell(y)$ strictly between 0 and 1 on all central summands. Ê(Recall that $\mathcal{M}$ has no abelian central summands.) ÊSo whenever $x \in \mathcal{M}, s_\ell(x) \perp s_\ell(y),$ \eqref{E:clarkson} gives
\begin{align*}
2(\|w B J(x) \|^p + \|T(y) \|^p) &= 2(\|(\begin{smallmatrix} x & 0 & \cdots
\end{smallmatrix})\|^p + \|y \|^p) \\
&= \|(\begin{smallmatrix} x Ê& 0 & \cdots \end{smallmatrix}) + y\|^p + \|(\begin{smallmatrix} x Ê& 0 & \cdots \end{smallmatrix}) - y\|^p \\
&= \|w B J(x) + T(y) \|^p + \|w B J(x) - T(y) \|^p,
\end{align*}
and this implies the \textit{right} supports of $J(x)$ and $T(y)$ are
orthogonal. ÊWrite $J = J_1 + J_2$ for the unique decomposition into multiplicative and antimultiplicative *-homomorphisms, with orthogonal ranges in $\mathcal{N}$. ÊSo
\begin{align*}
s_\ell(x) \leq s_\ell(y)^\perp &\Rightarrow s_r(T(y)) \perp [s_r(J_1(x)) + s_r(J_2(x))] \\ &\Rightarrow s_r(T(y)) \perp [J_1(s_r(x)) + J_2(s_\ell(x))].
\end{align*}
Since $s_r(x)$ can be any projection subequivalent to $s_\ell(y)^\perp$, and $J$ is unital, it follows in particular that $s_r(T(y)) \leq J_2(1)$. ÊThis inequality passes to the closed linear span of all $y$ under discussion, which is $R^p(\mathcal{M})(q-e_{11}^\mathcal{M})$ by an easy Hahn-Banach argument. ÊWe obtain
\begin{align*}
q'C^p(\mathcal{N}) &= T(R^p(\mathcal{M})q) \\ &= T(R^p(\mathcal{M})e_{11}^\mathcal{M}) + T(R^p(\mathcal{M})(q-e_{11}^\mathcal{M})) \\ &\subset \overline{wBJ(\mathcal{M})} + q'C^p(\mathcal{N})J_2(1) \\ & \subset s_\ell(w)C^p(\mathcal{N}) + q'C^p(\mathcal{N})J_2(1).
\end{align*}
Now $s_\ell(w)$ is equivalent to $e_{11}^\mathcal{N}$, so $(q' - s_\ell(w))$ has full central support. ÊMultiplying the containment above by $(q' - s_\ell(w))$ on the left, we get
$$(q' - s_\ell(w))C^p(\mathcal{N}) \subset (q' - s_\ell(w))C^p(\mathcal{N}) J_2(1).$$
This is only possible if $J_2(1) = 1_\mathcal{N}$, and so $J = J_2$ is antimultiplicative.

Finally, let $v$ be any partial isometry between orthogonal projections in $\mathcal{M}$. ÊWith $\pi: \mathcal{N} \hookrightarrow qM_\infty(\mathcal{M})q$ as before, we have
$$T( (\begin{smallmatrix} v^*v & 0 & \cdots \end{smallmatrix})) = w B J(v^*v) =
(w B J (v)) J(v^*) = T( (\begin{smallmatrix} v & 0 & \cdots
\end{smallmatrix})(\pi(J(v^*))).$$
Since $T$ is one-to-one, $(\begin{smallmatrix} v^*v & 0 & \cdots \end{smallmatrix}) = (\begin{smallmatrix} v & 0 & \cdots
\end{smallmatrix})(\pi(J(v^*)))$. ÊBut look at the left supports in $\mathcal{M}$ of these vectors; the first is $v^*v$ and the second is $\le s_\ell(v) = v v^*$. ÊThis contradiction finishes the proof.

\end{proof}

The only case remaining is an abelian central summand of $\mathcal{M}$ and
$\mathcal{N}$. ÊBecause column and row sums of $L^p(\mathbb{C}) = \mathbb{C}$
are identical, we cannot control the sizes of the commutants.

\begin{proposition}
Let $\mathcal{A} = L^\infty(X, \mu)$ be an abelian von Neumann algebra. ÊThe
$L^p$ $\mathcal{A}-\mathcal{A}$ bimodules are exactly the $p$-direct integrals
of measurable fields of Hilbert spaces over $(X, \mu)$.

By a $p$-direct integral we mean exactly the same construction as a direct
integral of a measurable field of Hilbert spaces, except that the norm is
$$\|\xi(\cdot)\| = \left( \int \|\xi(\omega)\|^p d\mu(\omega) \right)^{1/p}.$$
\end{proposition}

\begin{proof}
By Proposition \ref{T:center}, we may identify the left and right actions of $\mathcal{A}$. Ê(If presentations of $\mathcal{A}$ are given, this may involve an algebraic isomorphism.) ÊSubject to this, we claim that a right $L^p$ $\mathcal{A}$-module admits a \textit{unique} structure as an $L^p$ $\mathcal{A}-\mathcal{A}$ bimodule: left action of $\mathcal{A}$ given by $f \cdot \xi = \xi f$ and left inner product by $<\xi, \eta>_L = <\eta, \xi>_R.$

For suppose we are given an $L^p$ $\mathcal{A}-\mathcal{A}$ bimodule. ÊThat $f \cdot \xi = \xi f$ is automatic from the assumption; we further have, for any measurable set $E \subset X$,
\begin{align} \label{E:lr}
\| \chi_E <\xi, \xi>_L \| &= \|<\chi_E \xi, \chi_E \xi>_L \| = \|\chi_E \xi \|^2 \\ \notag &= \|\xi \chi_E \|^2 = \|<\xi \chi_E, \xi \chi_E>_R \| = \| <\xi, \xi>_R \chi_E \|.
\end{align}
Both $<\xi, \xi>_L$ and $< \xi, \xi >_R$ are positive functions in $L^{p/2}(X, \mu)$. ÊTaking $E$ in \eqref{E:lr} to be the set where one dominates the other, we deduce that $<\xi, \xi>_L $ $=$ $<\xi, \xi >_R$ $\mu$-a.e.. ÊBy polarization,
$$4 <\xi, \eta>_L = \sum_{k=1}^4 i^k <\xi + i^k \eta, \xi + i^k \eta>_L = \sum_{k=1}^4 i^k <\xi + i^k \eta, \xi + i^k \eta>_R = 4<\eta, \xi>_R.$$

So it is the same problem to describe the right $L^p$ $\mathcal{A}$-modules. ÊAny random projection $q: X \to \mathcal{P}(\mathcal{B}(\mathfrak{H}))$ gives a bimodule of the form $qC^p$, and since the isomorphism class of $qC^p$ depends only on the Murray-von Neumann equivalence class of $q$, we may assume that
$$q = \sum p_n \otimes \chi_{X_n} \in M_\infty (\mathcal{A}), \qquad \left( p_n
= \sum_{k=1}^n e_{kk} \right),$$
where $X_n$ is $\{\omega \mid \text{Tr} (q(\omega)) = n \}$.

By direct calculation, we see that $C^p(\mathcal{A})$ is the Bochner space $L^p( \ell^2, X, \mu)$, and
$$(p_n \otimes \chi_{X_n}) C^p \simeq L^p(\ell^2_n, X_n, \mu),$$
where we still use $\mu$ to denote the restricted measure on $X_n$. ÊThis last is a constant measurable
field of Hilbert spaces with norm
$$\left( \int_{X_n} \|f(\omega)\|^p d\mu(\omega) \right)^{1/p}.$$
The full module $q C^p$ is a central/$\ell^p$ sum of these,
$$qC^p = (\sum p_n \otimes \chi_{X_n}) C^p = \bigoplus_{\ell^p} L^p(\ell^2_n, X_n, \mu),$$
which is exactly a $p$-direct integral of the measurable field of Hilbert spaces which has dimension $n$ over $X_n$. ÊIt is clear that any such $p$-direct integral can be obtained in this way, so we are done.
\end{proof}

We summarize the results in
\begin{theorem} \label{T:bimodule}
Let $\mathcal{M}$ and $\mathcal{N}$ be $\sigma$-finite algebras, let $1 < p<\infty$, $p\ne 2$, and let $\mathfrak{X}$ be an $\mathcal{M}-\mathcal{N}$ $L^p$-bimodule which is countably generated and faithful for each action.

The centers of $\mathcal{M}$ and $\mathcal{N}$ are isomorphic and act identically on $\mathfrak{X}.$ ÊLet $z$ be the largest central projection which is abelian for both $\mathcal{M}$ and $\mathcal{N}$. ÊThen the left action of $z^\perp \mathcal{M}$ and the right action of $z^\perp \mathcal{N}$ are commutants on $z^\perp \mathfrak{X}.$ ÊOn the other hand, $z\mathfrak{X}$ is isomorphic to a $p$-direct integral of a measurable field of Hilbert spaces over $(X, \mu)$, where $z\mathcal{M} = z\mathcal{N} \simeq L^\infty(X, \mu).$
\end{theorem}

This has an appealing consequence.

\begin{theorem}
Under the same assumptions as above, there exists an $\mathcal{M}-\mathcal{N}$ $L^p$-bimodule if and only if $\mathcal{M}$ and $\mathcal{N}$ are Morita equivalent.
\end{theorem}

\begin{proof}
$\mathcal{M}$ and $\mathcal{N}$ are Morita equivalent exactly when $\mathcal{N}
\simeq qM_\infty(\mathcal{M})q$ for some projection $q$ with central support
one. Ê(This fact, and many other fundamental ideas, may be found in Rieffel's
discussions of Morita equivalence [Ri1], [Ri2].) ÊWhen this happens, one may
take $\mathfrak{X} \simeq R^p(\mathcal{M})q$ and notice $\mathcal{N} \simeq
(\mathcal{M}')^{\text{op}}.$

If there is an $\mathcal{M}-\mathcal{N}$ $L^p$-bimodule $\mathfrak{X}$, then
$\mathfrak{X} \simeq R^p(\mathcal{M}) q$ and $\mathcal{N} \simeq
\pi(\mathcal{N}) \subset qM_\infty(\mathcal{M})q.$ ÊLet $z$ be as in Theorem
\ref{T:bimodule}. ÊWe have that $\pi(z^\perp\mathcal{N}) = z^\perp
qM_\infty(\mathcal{M})q$, but
$$\pi(z\mathcal{N}) = z \mathcal{M} = ze_{11}M_\infty(\mathcal{M}) e_{11} \subset z q M_\infty(\mathcal{M}) q.$$ ÊThen
$$\mathcal{N} \simeq (z^\perp q + z e_{11})M_\infty(\mathcal{M})(z^\perp q + z
e_{11}),$$
so $\mathcal{M}$ and $\mathcal{N}$ are Morita equivalent.
\end{proof}

In fact, an $\mathcal{M}-\mathcal{N}$ $L^p$-bimodule $\mathfrak{X}$ which does
not degenerate on its abelian component (so $zq$ is abelian, and $L(\mathcal{M})$ and $R(\mathcal{N})$ are commutants) implements an equivalence of representation categories just as in the Hilbert C*-module case. ÊHere the densities are nonzero, and one makes use of the generalized relative tensor product, with functorial equivalence given by
$$Left\:L^q Mod (\mathcal{N}) \to Left\:L^r Mod(\mathcal{M}): {}_\mathcal{N}
\mathfrak{Y} \longmapsto \left( {}_\mathcal{M} \mathfrak{X}_\mathcal{N} \right)
\otimes_{\mathcal{N},p,q,r} \left({}_\mathcal{N} \mathfrak{Y} \right).$$

To see that this is an isomorphism, we let ${}_\mathcal{N}\bar{\mathfrak{X}}_\mathcal{M}$ be the contragredient and note that ``$({}_\mathcal{N}\bar{\mathfrak{X}}_\mathcal{M}) \otimes_{\mathcal{M},p,r,q}$" is the inverse map. ÊBy associativity of the relative tensor product, it suffices to show that
$$\mathfrak{X} \otimes_{\mathcal{N},p,q,r} \bar{\mathfrak{X}} \simeq L^p(\mathcal{M}); \qquad \bar{\mathfrak{X}} \otimes_{\mathcal{M},p,r,q} \mathfrak{X} \simeq L^p(\mathcal{N}).$$
We verify the first, using Proposition \ref{T:colxrow}:
$$\mathfrak{X} \otimes_{\mathcal{N},p,q,r} \bar{\mathfrak{X}} \simeq qC^p(\mathcal{N}) \otimes_{\mathcal{N},p,q,r} R^p(\mathcal{N})q \simeq qL^p(M_\infty({N}))q \simeq L^p(\mathcal{M}).$$
The second follows by symmetry.

\bigskip

\textit{Acknowledgments.} ÊThe first named author was partially supported by NSF grant DMS 00-88928. Ê The second named author would like to thank Stanis\l aw Goldstein for some helpful comments and also express special gratitude to Masamichi Takesaki; this paper is a lineal descendant of the dissertation written under his direction.

\end{document}